\documentclass[12pt]{article}

\usepackage{times}
\usepackage[a4paper,text={128mm,185mm}, centering]{geometry}
\usepackage{changepage}
\usepackage{titlesec}

\titleformat{\section}[hang]%
{\bfseries\large}{\thesection.}{1ex}{}%

\titleformat{\subsection}[hang]%
{\bfseries}{\thesubsection}{1ex}{}%

\usepackage{amsmath}
\usepackage{fancyhdr}
\usepackage{amsthm, amssymb, bm}
\usepackage[all]{xy}

\title{\vskip 5pt  \bf  THE TOPOLOGY OF CRITICAL PROCESSES,  II   $\hspace{-100pt}$
\newline \rm (THE FUNDAMENTAL CATEGORY)}

\author{\itshape\bfseries {Marco GRANDIS}}
\date{}

\def \LL {\skp \setlength{\leftskip}{8pt}}  
\def \LB {\setlength{\leftskip}{0pt}}   
\let \lan \langle
\let \ran \rangle
\def \sst {\scriptstyle}

\def \c {\colon}
\def \cc {\, \colon \!}

\def \q {\qquad}
\def \qq {\qquad \qquad}

\def \all {\forall}

\def \bu {{\scriptscriptstyle\bullet}}
\def \skp {\medskip}
\def \ndt {\noindent}
\def \Ndt {\medskip  \noindent}
\def\shp {^{\sharp}}
\def \tilde {\raise.17ex\hbox{$\scriptstyle\mathtt{\sim}$}}   
\def \ul {\underline}

\def \adj {\dashv}
\def \uw {\!\, \raisebox{0.3ex} {$\uparrow$}}
\def \dw {\raisebox{-0.2ex}{${\sst{\downarrow}}$}}

\def \rrw  {\; \raisebox{0.7ex}{$\longrightarrow$} \hspace{-3.7ex} \raisebox{-0.4ex}{$\longrightarrow$} \;}

\def \lrw {\;\; \raisebox{0.4ex}{$\longleftarrow$} \hspace{-3.7ex}
 \raisebox{-0.4ex}{$\longrightarrow$} \;\;}

\def \rlr { \;\; \raisebox{0.9ex}{$\longrightarrow$} \hspace{-3.7ex} 
 \raisebox{0.2ex}{$\longleftarrow$}
 \hspace{-3.7ex} \raisebox{-0.5ex}{$\longrightarrow$} \;\;}

\def \ard {\ar@{-->}}
\def \arp {\ar@{.>}}
\def \are {\ar@{->>}}
\def \aru {\ar@{=}}
\def \arv {\ar@{}}   
\def \arl {\ar@{-}}    
\def \arld {\ar@{--}}    
\def \arlp {\ar@{..}}    

\def \ti {\! \times \!}
\def \jo {{\, {\scriptstyle{\vee}} \,}}
\def \me {{\, {\scriptstyle{\wedge}}\, }}
\def \ci {{\raise.3ex\hbox{$  \scriptscriptstyle\circ $}}}
\def \Pro {\raisebox{0.45ex}{${\mbox{\fontsize{10}{10}\selectfont\ensuremath{\prod}}}$}}
\def \Sum {\raisebox{0.45ex}{${\mbox{\fontsize{10}{10}\selectfont\ensuremath{\sum}}}$}}
\def \Cup {\raisebox{0.4ex}{${\mbox{\fontsize{8}{10}\selectfont\ensuremath{\bigcup}}}$}}
\def \setm {{\raise.4ex\hbox{$ \, \scriptscriptstyle{\setminus} \; $}}}
\def \sub {\subset}
\def\le{\leqslant}
\def\ge{\geqslant}
\def \precc {\preceq}   
\def \iso {\: \cong \:}
\def \eq {\! \sim \!}
\def \eqq {\! \sim \!_2 \,}

\def \and {\mbox{ and }}
\def \for {\mbox{ for \,}}

\def \id {{\rm id\,}}
\def \op {^{{\rm op}}}
\def \pf  {^{{\rm pf}}}
\def \bbf  {^{{\rm bf}}}
\def \Im {{\rm Im\,}}
\def\dd {\partial}
\def\ddm {\partial^-}
\def\ddp {\partial^+}

\def \al {\alpha}

\def \ka {\kappa}
\def \si {\sigma}
\def \ph {\varphi}

\def \two {{\bf 2}}
\def \three {{\bf 3}}
\def \N {{\bf N}}   
\def \Z {{\bf Z}}   
\def \R {{\bf R}}   
\def\bc {{\bf c}}

\def \Cat {\mathsf{Cat}}
\def \Top {\mathsf{Top}}

\def \dTop {{\rm d}\mathsf{Top}}
\def \cTop {{\rm c}\mathsf{Top}}
\def \cpfTop {{\rm c_{pf}}\mathsf{Top}}
\def \cbfTop {{\rm c_{bf}}\mathsf{Top}}

\def \Fl {\mathsf{F}{\rm l}\,}
\def \uPi {\uw\Pi}
\def \upi {\uw\pi}

\def \sing  {\{*\}}
\def \bbN {{\mathbb{N}}}  
\def \bbZ {\mathbb{Z}}
\def \bbI {{\mathbb{I}}}   
\def \bbR {{\mathbb{R}}}
\def \bbS {\mathbb{S}}

\def \uI {{\uw\mathbb{I}}}
\def \uR {\uw\mathbb{R}}
\def \uS {\uw\mathbb{S}}

\def \rc {{\rm c}}
\def \cI {{\rm c}\mathbb{I}}
\def \cJ {{\rm c}\mathbb{J}}
\def \cR {{\rm c}\mathbb{R}}
\def \cS {{\rm c}\mathbb{S}}

\begin{document}
\maketitle

\vskip 25pt
\begin{adjustwidth}{0.5cm}{0.5cm}
{\small
{\bf Abstract.} Directed Algebraic Topology studies spaces equipped with a form of direction, to include 
models of non-reversible processes. In the present extension we also want to cover {\em critical 
processes}, indecomposable and unstoppable.

	The first part of this series introduced {\em controlled spaces}, examining how they can model 
critical processes in various domains, from the change of state in a memory cell to the action of a 
thermostat or a siphon. We now construct the fundamental category of these spaces.

	[The previous version of this Preprint was too long and has been split in Parts II and III.]\\

\vskip -8pt \noindent
{\bf Keywords.} Directed algebraic topology, homotopy theory, fundamental category.\\

\vskip -8pt \noindent
{\bf Mathematics Subject Classification (2010).} 55M, 55P, 55Q.
}
\end{adjustwidth}


\section*{Introduction}\label{Intro}

\subsection{Critical processes and controlled spaces}\label{0.1}
Directed Algebraic Topology is an extension of Algebraic Topology, dealing with `spaces' where the paths 
need not be reversible; the general aim is including the representation of {\em irreversible processes}. A 
typical setting for this study, the category $\dTop$ of directed spaces, or d-spaces, was introduced and 
studied in \cite{G1}--\cite{G3}; it is frequently employed in the theory of concurrency: see the book 
\cite{FGHMR} and many articles cited in the previous paper \cite{G5}.

	The present series is devoted to a further extension, where the paths can also be non-decomposable 
in order to include {\em critical processes}, indivisible and unstoppable -- either reversible or not. 
For instance: quantum effects, the onset of a nerve impulse, the combustion of fuel in a piston, the switch 
of a thermostat, the change of state in a memory cell, the action of a siphon, moving in a no-stop road, etc.

	To this effect the category of d-spaces was extended in Part I \cite{G5} to the category $\cTop$ of 
{\em controlled spaces}, or {\em c-spaces}: an object is a topological space equipped with a set $X\shp$ 
of continuous mappings $a\c [0, 1] \to X$, called {\em controlled paths}, or {\em c-paths}, that satisfies 
three axioms:

\LL
	
\ndt (csp.0) ({\em constant paths}) the trivial loops at the endpoints of a controlled path are controlled,
	
\Ndt (csp.1) ({\em concatenation}) the concatenation of consecutive controlled paths is controlled,
	
\Ndt (csp.2) ({\em global reparametrisation}) the reparametrisation of a controlled path by a surjective 
increasing map $ [0, 1] \to [0, 1] $ is controlled.

\LB

\skp 	A {\em map of c-spaces}, or {\em c-map}, is a continuos mapping which preserves the selected paths. 
Their category $\cTop$ contains the category $\dTop$ of d-spaces as a full subcategory, reflective and 
coreflective: a c-space is a d-space if and only if it is {\em flexible}, which means that each point is flexible 
(its trivial loop is controlled) and every controlled path is flexible (all its restrictions are controlled).

	Here we deal with the fundamental category of controlled spaces. Part III will study more advanced 
methods of computations of the latter, with applications to models of critical processes and concurrency. 
The homotopy theory of c-spaces will be dealt with in Part IV.

\subsection{Standard intervals}\label{0.2}
The difference between these two settings -- directed and controlled spaces -- shows clearly in two structures 
of the euclidean interval $ \bbI = [0, 1]$, the starting point of homotopy in each setting, represented as follows 
(and better described in \ref{1.3})
%
    \begin{equation}
\xy <.5mm, 0mm>:
(0,-6) *{\sst{0}}; (40,-6) *{\sst{1}}; (110,0) *{\bu}; (150,0) *{\bu}; (60,-3) *{\uI}; 
(110,-6) *{\sst{0}}; (150,-6) *{\sst{1}}; (170,-3) *{\cI}; 
(0,5) *{}; (0,-10) *{};
\POS(0,0) \arl+(40,0), \POS(19,0) \ar+(4,0),
\POS(110,0) \arl+(40,0), \POS(0,-1) \arl+(0,2), \POS(40,-1) \arl+(0,2), 
\POS(129,0) \are+(4,0),
\endxy
    \label{0.2.1} \end{equation}

	The left figure shows the {\em standard d-interval} $\uI$ (in $\dTop$): its directed paths are all the 
(weakly) increasing maps $ \bbI \to \bbI$. It may be viewed as the essential model of a non-reversible 
process, or a one-way route in transport networks. Its fundamental category $ \uPi_1(\uI)$, as defined in 
\cite{G1, G3}, is the ordered set $ [0, 1]$, with one arrow $ t \to t' $ for each pair $ t \le t' $ in $ [0, 1]$.

	The right figure shows the {\em standard c-interval} $ \cI$, or {\em one-jump interval} (in $ \cTop$): 
its controlled paths are the surjective increasing maps $ \bbI \to \bbI $ and the trivial loops at 0 or 1. It models 
a non-reversible unstoppable process, or a one-way no-stop route. Its fundamental category $ \uPi_1(\cI)$, 
as defined here, is the ordinal $\two$, with one non-trivial arrow $ 0 \to 1 $ and two identities.

\subsection{Outline}\label{0.3}
The basic definitions and the main examples of Part I are recalled in Section 1; in particular, every c-space 
has two associated d-spaces, the generated d-space $\hat{X}$ and the flexible part $\Fl X$, by the reflector 
and coreflector of the embedding $\dTop \to \cTop$ (see \ref{1.2}). Then we introduce in
Section 2 two weak forms of flexibility which will play a role in our study: preflexible and border flexible 
c-spaces.

	Section 3 reviews the basic part of the homotopy theory of d-spaces studied in \cite{G1, G3}, including 
the construction of the fundamental category $\uPi_1(X)$ of a d-space; some of these results are 
already extended or adapted to c-spaces.

	In the next two sections we introduce the fundamental category $\uPi_1(X)$ of a c-space, an extension 
of the previous case. Its vertices are the flexible points of $ X$; its arrows come out of a complex construction, 
based on the hybrid square $\cI \ti \uI$: they are equivalence classes of controlled paths (parametrised on
$\cI$), up to flexible deformations (parametrised on $\uI$). 

	The main results can be found in Section 5: the 
construction of $\uPi_1(X)$, its relationship with the fundamental category of the associated d-spaces 
$\hat{X}$ and $\Fl X$ (in  \ref{5.2} and  \ref{5.3}), its homotopy invariance (in \ref{5.4}) 
and its calculation for covering maps of c-spaces (in \ref{5.8}). All this is based on the technical 
analysis of Section 4. We end by calculating the fundamental category of the basic c-spaces, in \ref{5.9}.

\subsection{Comments}\label{0.4}
In the present extension we reach models of phenomena that have no place in the previous settings of 
Directed Algebraic Topology, and peculiar `shapes', like the {\em one-stop circle} $ \cS^1$, the 
{\em $n$-stop circle} $ \rc_n\bbS^1$, or the higher controlled spheres and tori described in Part I. The 
fundamental category of the new spaces is often quite simple.

	We also loose some good properties of the theory of d-spaces. For instance, the interval $\cI$ is not 
exponentiable in $ \cTop $ (see \ref{4.7}(b)), and the associated cylinder functor $\, I(X) = X \ti \cI \,$ has no 
right adjoint: there is no path endofunctor. We manage to extend the fundamental category, by allowing 
c-paths to be deformed by flexible homotopies, and one can also extend directed singular homology, but new 
methods of calculation are needed: the van Kampen theorem and the Mayer-Vietoris sequence are based 
on the subdivision of paths and homological chains, which is no longer permitted. On the other hand, the 
theory of covering maps can be extended to the present case (Theorem \ref{5.8}(b)). 

	Essentially, the previous setting of d-spaces extends classical topology by breaking the symmetry of 
reversion: directed paths need no longer be reversible. This further extension to c-spaces breaks a flexibility 
feature that d-spaces still retain: paths can no longer be subdivided, and this has drastic consequences.

	Critical processes and transport networks are often represented by graphs, in an effective way as far 
as they do not interact with continuous variation. We want to show that they can also be modelled by 
structured spaces, in a theory that includes classical topology and `non-reversible spaces'. Controlled 
spaces can thus unify aspects of continuous and discrete mathematics, interacting with hybrid control systems 
and others sectors of Control Theory \cite{Br, He}.

\subsection{Notation and conventions}\label{0.5}
A continuous mapping between topological spaces is called a {\em map}. $ \bbR $ denotes the euclidean 
line as a topological space, and $ \bbI $ the standard euclidean interval $ [0, 1]$. The identity path $ \id \bbI $ 
is written as $ \ul{i}$. The open and semiopen intervals of the real line are denoted by square brackets, 
like $ ]0, 1[$, $[0, 1[ $ etc. A space is {\em locally compact} if every point has a local basis of compact 
neighbourhoods; the Hausdorff property is not assumed.

	A {\em preorder} relation, generally written as $ x \prec y$, is assumed to be reflexive and transitive; 
an {\em order}, often written as $ x \le y$, is also assumed to be anti-symmetric. A mapping which preserves 
(resp.\ reverses) preorders is said to be {\em increasing} (resp.\ {\em decreasing}), always used in the weak 
sense. 

	As usual, a preordered set $ X $ is identified with the small category whose objects are the elements 
of $ X$, with one arrow $ x \to x' $ when $ x \prec x' $ and none otherwise.
 
	The binary variable $\al$ takes values $ 0, 1, $ which are generally written as $ -, + $ in superscripts 
and subscripts. The symbol $\, \sub \,$ denotes weak inclusion. 

	The first paper \cite{G5} of this series is cited as Part I; the reference I.2 or I.2.3 points to Section 2 
or Subsection 2.3 of Part I, respectively.

\section{Controlled and directed spaces}\label{s1}

	We recall the definition of controlled space, or c-space, introduced in Part I. Their category $ \cTop $ is 
an extension of the category $\dTop$ of directed spaces, or d-spaces, studied in \cite{G1}--\cite{G3} and 
commonly used in concurrency; we generally refer to the book \cite{G3}. The term `selected path' is used 
in both cases.

\subsection{Controlled spaces}\label{1.1}
	As defined in Part I, a {\em controlled space} $ X$, or {\em c-space}, is a topological space equipped 
with a set $ X\shp $ of (continuous) maps $ a\c \bbI \to X$, called {\em controlled paths}, or {\em c-paths}, 
that satisfies three axioms:

\LL

\ndt (csp.0) ({\em constant paths}) the trivial loops at the endpoints of a controlled path are controlled,

\Ndt (csp.1) ({\em concatenation}) the controlled paths are closed under path concatenation: if the 
consecutive paths $ a, b $ are controlled, their concatenation $ a * b $ is also,

\Ndt (csp.2) ({\em global reparametrisation}) the controlled paths are closed under pre-composition with 
every surjective increasing map $ \rho\c \bbI \to \bbI$: if $ a $ is a controlled path, $a\rho$ is also.

\LB

\skp 	The controlled paths are also closed under general $n$-ary concatenations, based on arbitrary 
partitions $ 0 = t_0 < t_1 < ... < t_n = 1 $ of $ \bbI $ (as shown in I.1.2). The underlying topological space 
is written as $U(X)$, or $|X|$, and called the {\em support} of $X$.

	A {\em map of c-spaces}, or {\em c-map}, is a continuos mapping which preserves the selected paths. 
Their category is written as $\cTop$.

	Reversing c-paths, by the involution $r(t) = 1 - t$, yields the {\em opposite} c-space $RX = X\op$, where 
$a \in (X\op)\shp$ if and only if $ar$ belongs to $ X\shp$. We have thus the {\em reversor} endofunctor
    \begin{equation}
R\c \cTop \to \cTop,   \qq   RX  =  X\op.
    \label{1.1.1} \end{equation}

	The c-space $ X $ is {\em reversible} if $ X = X\op$. More generally, it is {\em reversive} if it is 
isomorphic to $ X\op$.

	Controlled spaces have all limits and colimits, which are topological limits and colimits with the initial 
or final structure determined by the structural maps.

	In a c-space $X$, a point $ x $ is {\em flexible} if its trivial loop $ e_x $ is controlled. The {\em flexible 
support} $ |X|_0 $ is the subspace of these points, and can be empty. A c-path is {\em flexible} if all its 
restrictions are controlled. The c-space itself is {\em flexible} if every point and every selected path are 
flexible.

	The singleton space has two c-structures: the discrete one $D_c\sing$, with no controlled paths, and the 
{\em flexible singleton} $ \sing $ where the trivial loop is selected; this is the terminal object and the unit of 
the cartesian product. 

	The category $ \cTop $ has two forgetful functors to topological spaces
    \begin{equation}
U\c \cTop \to \Top,   \qq   U_0\c \cTop \to \Top,
    \label{1.1.2} \end{equation}
where $ U(X) = |X| $ is the topological support and $ U_0(X) = |X|_0 $ is the flexible support. $ U $ has 
both adjoints, $ U_0 $ has only the left one: it preserves limits and sums (I.1.7(e)).

	We say that the c-space $X_1$ is {\em finer} than $ X_2 $ if they have the same topological support 
$ X $ and $X_1\shp \sub X_2\shp$, which means that the identity map of $ X $ is a c-map $ X_1 \to X_2$; 
the latter is called a {\em reshaping}.

	We shall see later, in \ref{4.9}, that requiring that {\em all trivial loops be controlled} would be a serious 
hindrance.

\subsection{Directed spaces}\label{1.2}
Previously, a {\em directed space} $ X$, or {\em d-space}, was also defined as a topological space with a 
set $ X\shp $ of selected paths, called {\em directed paths}, or {\em d-paths}, under stronger axioms: the 
selected paths are stable under: (all) constant paths, concatenations and {\em partial} reparametrisations 
(by increasing endomaps of the interval, not assumed to be surjective) \cite{G3}.

	A d-space is the same as a flexible c-space, and can also be defined in this way. The category 
$\dTop$ of d-spaces and d-maps is a full subcategory of $ \cTop$, reflective and coreflective
    \begin{equation} \begin{array}{ccr} 
\hat{\, }\c  \cTop \to \dTop   &&    \mbox{(the {\em reflector}),}
\\[5pt]
\Fl\c  \cTop \to \dTop   &\q&    \mbox{(the {\em coreflector}).}
    \label{1.2.1} \end{array} \end{equation}

	In the first construction the {\em generated d-space} $\hat{X}$ of a c-space $ X $ has the same support 
with the d-structure generated by the c-paths, i.e.\ the finest containing them; it can be obtained by 
stabilising the latter under constant paths, restriction and general concatenation. The unit of the adjunction 
is the reshaping $ X \to \hat{X}$; the counit is the identity $\hat{Y} = Y$, for a d-space $ Y$.

	In the second construction the {\em flexible part} $\Fl X $ is the flexible support $ |X|_0 $ with the 
d-structure of the flexible c-paths. The counit is the inclusion $\Fl X \to X$, the unit is the identity 
$Y = \Fl Y$, for a d-space $Y$.

	Also $\dTop$ has all limits and colimits, preserved by the embedding in $ \cTop$.

\subsection{Structured intervals and cubes}\label{1.3}
(a) In $\dTop$ the {\em standard d-interval} $\uI$ has the d-structure generated by the identity $\ul{i} = \id \bbI$: 
the directed paths are all the increasing maps $\bbI \to \bbI$. It plays the role of the standard interval in 
this category, because the directed paths of any d-space $X$ coincide with the d-maps $\uI \to X$.

	It may be viewed as an essential model of a non-reversible process, or a one-way road in transport 
networks. It will be represented as in figure \eqref{0.2.1} of the Introduction.

	Similarly, the directed line $\uR$ has for directed paths all the increasing maps $\bbI\to \bbR$.

\Ndt (b) In $ \cTop $ the {\em standard c-interval} $ \cI$, or {\em one-jump interval}, has the same support, 
with the c-structure generated by the identity $ \ul{i}$: the controlled paths are the surjective increasing 
maps $ \bbI \to \bbI $ and the trivial loops at 0 or 1. The controlled paths of any c-space $ X $ coincide with 
the c-maps $ \uI \to X$.

	It can model a {\em non-reversible unstoppable process}, or a {\em one-way no-stop road}. It is also 
represented in figure \eqref{0.2.1}, marking by a bullet the isolated flexible points: in this case, the endpoints 
of the interval.

	The controlled line $ \cR $ has for directed paths all the increasing maps $ \bbI \to \bbR $ whose image 
is an interval $ [k, k'] $ with integral endpoints.

	We shall also use other c-structures of the compact interval, already examined in I.2.4.

\Ndt (c) In the {\em two-jump interval} $ \cJ$, the non-trivial c-paths are the increasing maps $ \bbI \to \bbI $ 
whose image is either $ [0, 1/2]$, or $ [1/2, 1]$, or $ [0, 1]$
%
    \begin{equation}
\xy <.5mm, 0mm>:
(0,-6) *{\sst{0}}; (40,-6) *{\sst{1/2}}; (80,-6) *{\sst{1}}; 
(0,0) *{\bu}; (40,0) *{\bu}; (80,0) *{\bu}; (100,-3) *{\cJ}; 
(0,5) *{}; (0,-10) *{};
\POS(0,0) \arl+(80,0), \POS(20,0) \are+(4,0), \POS(60,0) \are+(4,0),
\endxy
    \label{1.3.1} \end{equation}

	This c-space can model a non-reversible two-stage process. Formally, it parametrises the ordinary 
concatenation of two c-paths, see \ref{4.2}.

\Ndt (d) The {\em reversible c-interval} $ \cI^\sim$ has a c-structure generated by the identity $ \ul{i} $ and 
the reversion $ r\c \bbI \to \bbI$. It can model a reversible unstoppable process. The reversible c-paths of a 
c-space $ X $ coincide with the c-maps $ \cI^\sim \to X$.

\Ndt (e) We shall also use the {\em delayed intervals} $\rc_-\bbI $ and $\rc_+\bbI$ of I.2.4(b). Each of these 
c-structure is generated by a single map $ \bbI \to \bbI$, namely $ \rho(t) = 0 \jo (2t - 1) $ or 
$ \si(t) = 2t \me 1$, respectively. They can model irreversible non-stoppable processes with inertia, or 
inductance.

\Ndt (f) The powers $ \bbI^n $ and $ \bbR^n $ inherit various controlled structures. On $ \bbI^2 $ we 
shall mostly use the d-square $ \uI^2$, the c-square $ \cI^2 $ and the hybrid square $ \cI \ti \uI$ (cf.\ I.2.7).

\subsection{Structured spheres}\label{1.4}
(a) The {\em standard d-circle} $ \uS^1 $ can be obtained as an orbit space
%
    \begin{equation} 
\xy <1pt, 0pt>:
(100,-12) *{\uS^1 = (\uR)/\bbZ}; (0,30) *{}; (0,-30) *{};
%
(0,0) = "A"  *\cir<25pt>{}="ca"  
\POS(-1,25)\ar-(2,0)
\endxy
    \label{1.4.1} \end{equation}
with respect to the action of the group of integers on the directed line $ \uR $ (by translations): the 
directed paths of $ \uS^1 $ are the projections of the increasing paths in the line.

	$\uS^1 $ can also be obtained as the coequaliser in $\dTop$ of the following pair of maps
    \begin{equation}
\ddm, \ddp\c  \sing  \rrw  \uI,   \qq    \ddm(*)  =  0,   \;\;\;   \ddp(*)  =  1,
    \label{1.4.2} \end{equation}
that is the quotient $ \uI/\dd \bbI$ which identifies the points of the boundary $\dd \bbI = \{0, 1\}$.

\Ndt (b) More generally, the {\em directed $n$-dimensional sphere} is defined, for $ n > 0$, as the quotient 
of the directed cube $ \uI^n $ modulo the equivalence relation that collapses its boundary $ \dd \bbI^n $ 
to a point
    \begin{equation}
\uS^n  =  (\uI^n)/(\dd \bbI^n) \;\;\;  (n > 0),   \q   \uS^0  =  \bbS^0  =  \{-1, 1\},
    \label{1.4.3} \end{equation}
while $ \uS^0 $ has the discrete topology and the natural (discrete) d-structure.

\Ndt (c) The {\em standard c-circle} $ \cS^1$, or {\em one-stop circle}, can also be defined as an orbit space
%
    \begin{equation} 
\xy <1pt, 0pt>:
(25,0) *{\bu}; (100,-12) *{\cS^1 =  (\cR)/\bbZ}; (0,30) *{}; (0,-30) *{};
%
(0,0) = "A"  *\cir<25pt>{}="ca"  
\POS(-1,24.8)\are-(4,0)
\endxy
    \label{1.4.4} \end{equation}
for the action of the group of integers on the controlled line $ \cR$, by translations. The controlled paths 
of $ \cS^1 $ are the projections of the controlled paths in the line: here this means an anticlockwise path 
(in the ordered plane) which is a loop at $* = [0]$, the only flexible point (corresponding to $ (1, 0) $ in the 
plane).

	The c-space $ \cS^1 $ can also be obtained as the coequaliser in $ \cTop $ of the following pair of 
maps 
    \begin{equation}
\ddm, \ddp\c \sing  \rrw \, \cI,   \q   \ddm(*)  =  0,   \;\;\;   \ddp(*)  =  1.
    \label{1.4.5} \end{equation}

\Ndt (d) More generally, the {\em $n$-stop c-circle} $\rc_n\bbS^1$ $(n > 0)$ is constructed in I.2.6(b) as the 
orbit space
    \begin{equation}
\rc_n\bbS^1  =  (\rc_n\bbR)/\bbZ   \qq   (\rc_1\bbS^1  =  \cS^1),
    \label{1.4.6} \end{equation}
where the c-paths of $ \rc_n\bbR $ are the increasing paths whose image is an interval $ [k/n, k'/n]$, for 
integers $ k \le k'$.

\Ndt (e) The {\em standard c-sphere} $ \cS^n $ is defined as a quotient of the cube $\cI^n$ (for $n > 0$)
    \begin{equation}
\cS^n  =  (\cI^n)/(\dd \bbI^n)   \;\;\;  (n > 0),   \q   \cS^0  =  \bbS^0  =  \{-1, 1\}.
    \label{1.4.7} \end{equation}
%
%

\subsection{Identities and associativity}\label{1.5}
Concatenation of paths and the various forms of reparametrisation have a complex relationship. 
Here we recall two well-known points.  

\Ndt (a) The constant loops act as identities up to the equivalence relation generated by global 
reparametrisation. In fact, the following surjective increasing maps $ \rho, \si\c \bbI \to \bbI $ 
reparametrise any path $ a$, from $ $x to $ y$, as $ e_x * a $ or $ a * e_y$, respectively:
%
    \begin{equation} 
\xy <.5mm, 0mm>:
%
(-10,20) *{}; (215,-20) *{};
(40,-10) *{\rho}; (110,-10) *{\si}; 
(165,10) *{\rho(t)  =  0 \jo (2t - 1),}; (155,-5) *{\si(t)  =  2t \me 1,}; 
@i@={(0,-15), (30,-15), (30,15), (0,15)},
s0="prev"  @@{;"prev";**@{-}="prev"};
@i@={(70,-15), (100,-15), (100,15), (70,15)},
s0="prev"  @@{;"prev";**@{-}="prev"};
\POS(0,-14.5) \arl+(15,0),  \POS(15,-15) \arl+(15,30),
\POS(70,-14.5) \arl+(15,30),  \POS(85,14.5) \arl+(15,0),
\endxy
    \label{1.5.1} \end{equation}
    \begin{equation*}
a\rho  =  e_x * a,    \qq    a\si  =  a * e_y.
    \label{} \end{equation*}

	These maps were used in I.2.4(b), as the {\em past-} (resp.\ {\em future-}) 
{\em delayed reparametrisation}. Since $ \rho \le \ul{i} \le \si$, there are directed homotopies with fixed 
endpoints $a\rho \to a \to a\si$ in $\dTop$ \cite{G3}, which also work in the present setting, as we shall 
see in Lemma \ref{4.6}. (These homotopies are reversible in $ \Top$, but are not for our directed structures.)

\Ndt (b) All $n$-ary concatenations are equivalent, up to invertible reparametrisation (cf.\ I.1.2). In 
particular, let $ \rho\c \bbI \to \bbI $ be the obvious piecewise affine invertible reparametrisation that takes 
the partition $ 0 < 1/2 < 3/4 < 1 $ to the regular partition $ 0 < 1/3 < 2/3 < 1$, while $ \si\c \bbI \to \bbI $ does 
the same on $ 0 < 1/4 < 1/2 < 1$
%
    \begin{equation} 
\xy <.5mm, 0mm>:
%
(0,20) *{}; (150,-20) *{};
(15,-24) *{\sst{1/2}}; (-10,-5) *{\sst{1/3}}; (-10,5) *{\sst{2/3}}; 
(115,-24) *{\sst{1/2}}; (90,-5) *{\sst{1/3}}; (90,5) *{\sst{2/3}}; 
(42,-10) *{\rho}; (142,-10) *{\si}; 
@i@={(0,-15), (30,-15), (30,15), (0,15)},
s0="prev"  @@{;"prev";**@{-}="prev"};
@i@={(100,-15), (130,-15), (130,15), (100,15)},
s0="prev"  @@{;"prev";**@{-}="prev"};
\POS(0,-15) \arl+(15,10),  \POS(15,-5) \arl+(15,20),
\POS(100,-15) \arl+(15,20),  \POS(115,5) \arl+(15,10),
\POS(15,-16) \arl+(0,2), \POS(22.5,-16) \arl+(0,2), 
\POS(-1,-5) \arl+(2,0), \POS(-1,5) \arl+(2,0), 
\POS(107.5,-16) \arl+(0,2), \POS(115,-16) \arl+(0,2), 
\POS(99,-5) \arl+(2,0), \POS(99,5) \arl+(2,0), 
\endxy
    \label{1.5.2} \end{equation}

	Now, if $ d = a * b * c $ is the regular concatenation of three consecutive paths, based on the 
	partition $ 0 < 1/3 < 2/3 < 1$
    \begin{equation}
d\rho  =  a * (b * c),    \qq    d\si  =  (a * b) * c.
    \label{1.5.3} \end{equation}

	Again $ \rho \le \ul{i} \le \si$, and there are homotopies with fixed endpoints $ d\rho \to d \to d\si$, 
which work in $\dTop$ and will also work in the present setting, by Lemma \ref{4.6}.

\section{Weak flexibility}\label{s2}

	We now introduce weak forms of flexibility that will be important for the construction of the 
fundamental category, and still hold in basic c-spaces like $ \cI$, $\cJ$, $\cR$, $\cS^1$ and all their 
products (and limits) -- although all of them are rigid c-spaces, in the sense of I.1.6.

	$X $ is always a c-space.

\subsection{Preflexible and border flexible c-spaces}\label{2.1}
(a) If $ S $ is a subset of the flexible support $ |X|_0 $ of $ X$, we can form a finer c-space $ X_{|S} $ on 
the same support selecting the c-paths of $ X $ whose endpoints belong to $ S$. We say that the c-space 
$ X_{|S} $ is {\em full} in $ X$, or {\em obtained from} $ X $ {\em by restricting the flexible support to} $ S$. 
We shall see, in Theorem \ref{5.3}(a), that $ \uPi_1(X_{|S}) $ is the full subcategory of $ \uPi_1(X) $ 
with vertices in $ S$. 

	For instance the c-spaces $ \cI $ and $ \cJ $ are full in $\uI$, which is just finer than $ \bbI$. 
Both relationships, being finer or full, are preserved by products.

\Ndt (b) We say that $ X $ is {\em preflexible} if it is full in the generated d-space $\hat{X}$, which means that 
every c-path of $ \hat{X} $ between flexible points of $ X $ is already a c-path of the latter. The interest of 
this property will be evident in Theorem \ref{5.3}(b).

\Ndt (c) We say that $ X $ is {\em border flexible} if its controlled paths are closed under `cutting out delays at 
the endpoints'. More precisely, we are asking that, for every c-path $ a $ of $ X $ which is constant on 
$ [0, t_1] $ and $ [t_2, 1]$, the {\em border restriction} $a\rho$ be still a controlled path, for
    \begin{equation}
\rho\c \bbI \to \bbI,   \q   \rho(t)  =  t_1 + (t_2 - t_1)t   \q\;\;   (0 \le t_1 < t_2 \le 1).
    \label{2.1.1} \end{equation}

	Plainly, every preflexible c-space is border flexible; the converse is false: see \ref{2.2}(b).

\Ndt (d) We also introduce the {\em path-support} $ |X|_1 $ of the c-space $ X $ as the topological subspace 
of $ |X| $ formed by the union of the images of all c-paths in $ X$, so that $ |X|_0 \sub |X|_1 \sub |X|$. A 
c-map can be restricted to the path-supports. We say that $ X $ has a {\em total path-support} if $|X|_1 = |X|$.

\subsection{Examples and remarks}\label{2.2}
(a) Besides all d-spaces, many c-spaces we have considered in Part I and here are preflexible: for instance 
$ \cI$, $\cJ $ and $ \cI^\sim $ (in $ \uI$), $ \cR $ (in $ \uR$), $ \cS^1 $ (in $ \uS^1$), and all their limits 
(by Proposition \ref{2.3}).

	The delayed intervals $\rc_-\bbI $ and $\rc_+\bbI $ recalled in \ref{1.3}(e) are not even border flexible: 
the preflexible space generated by any of them (according to Proposition \ref{2.3}(a)) is the standard interval 
$ \cI$.

	It is not difficult to prove that the c-spheres $ \cS^n $ are not border flexible, for $n \ge 2$.

\Ndt (b) The `diagonal' c-structure $ X $ of the square $ \bbI^2 $ described in I.2.7(d) is border flexible 
(obviously) and not preflexible (as shown below); it does not have a total path-support.

	We recall that the c-paths of $ X $ are generated by two diagonal paths, $ t \mapsto (t, t) $ and 
$ t \mapsto (t, 1 - t)$, represented in the left figure below
%
    
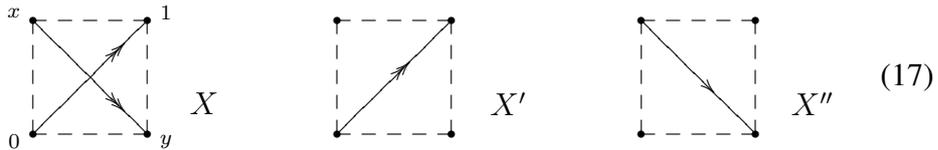
\begin{equation} 
\xy <.5mm, 0mm>:
%
(0,22) *{}; (0,-22) *{};
(45,-7) *{X}; (125,-7) *{X'}; (205,-7) *{X''}; 
(-5,-17) *{\sst{0}};  (-5,17) *{\sst{x}}; (35,-17) *{\sst{y}}; (35,17) *{\sst{1}};
@i @={(0,-15), (30,-15), (30,15), (0,15),
(80,-15), (110,-15), (110,15), (80,15), 
(160,-15), (190,-15), (190,15), (160,15)} @@{*{\bu}};
@i @={(0,-15), (30,-15), (30,15), (0,15)},
s0="prev"  @@{;"prev";**@{--}="prev"};
@i @={(80,-15), (110,-15), (110,15), (80,15)},
s0="prev"  @@{;"prev";**@{--}="prev"};
@i @={(160,-15), (190,-15), (190,15), (160,15)},
s0="prev"  @@{;"prev";**@{--}="prev"};
\POS(0,-15) \arl+(30, 30), \POS(0,15) \arl+(30, -30), 
\POS(80,-15) \arl+(30, 30), \POS(160,15) \arl+(30, -30), 
\POS(22,7) \are+(2,2), \POS(22,-7) \are+(2,-2),
\POS(97,2) \are+(2,2), \POS(177,-2) \ar+(2,-2), 
\endxy
    \label{2.2.1} \end{equation}

	The flexible points are the four vertices of the square. The generated d-space $ \hat{X} $ has also 
c-paths $ 0 \to y $ and $ x \to 1$, proving that $ X $ is not preflexible.

	We also note that the structure of $ X $ is generated by two finer c-structures $ X'$, $X'' $ of the square, 
with the same flexible points and non-trivial c-paths generated by one of the previous diagonals. $ X $ is 
thus the pushout of three preflexible spaces, $ X' $ and $ X'' $ over $ X_0$: the latter is the square 
$ \bbI^2 $ with the trivial loops at the vertices (the intersection of the structures of $ X' $ and $ X''$).

\Ndt (c) The quotient of the interval $ X = \rc[0, 2] \sub \cR $ modulo the equivalence relation that collapses 
its second half to a point
    \begin{equation}
X/[1, 2]  \iso  \rc_+\bbI,
    \label{2.2.2} \end{equation}
is not border flexible. It is the pushout in $ \cTop $ of two maps of border flexible c-spaces, the inclusion 
$\rc[1, 2] \to \rc[0, 2] $ and the map $ \rc[1, 2] \to \sing$.

\Ndt (d) {\em Remarks}. In a border flexible c-space, initial or final delays cannot be required; but let us note 
that they can never be prevented -- a global reparametrisation can always produce them.

	On the other hand, `internal' delays can be required, as in the border-flexible {\em middle-delay} interval 
$\rc_M\bbI$, with the c-structure generated by the map 
$$ \rc_M(t)  =  ((3t \me 1) \jo (3t - 1))/2 \qq    \mbox{({\em middle-delay map}).}$$

\ndt (e) Full c-spaces were already considered in I.3.5, in the equivalent perspective of excluding the 
flexible points of $|X|_0 \setm S$.

\subsection{Proposition and Definition {\rm (Reflectors and limits)}}\label{2.3}
{\em
(a) Preflexible c-spaces form a full, reflective subcategory $\cpfTop$ of $\cTop$. The reflector
    \begin{equation}
( - )\pf\c  \cTop \to \cpfTop,   \qq   X  \mapsto  X\pf,
    \label{2.3.1} \end{equation}
gives the {\em generated preflexible space} $ X\pf$, with the same flexible points and all the c-paths of 
$ \hat{X} $ between them; maps stay `unchanged'. The unit is the reshaping $ X \to X\pf$.

	As in every full reflective subcategory, preflexible c-spaces are closed in $ \cTop $ under limits. 
All colimits can be obtained taking the colimit in $ \cTop $ and applying the reflector. Sums are preserved 
by the inclusion, but pushouts are not -- in general.

\Ndt (b) Similarly, border flexible c-spaces form a full, reflective subcategory $\cbfTop$ of $ \cTop$. 
The reflector
    \begin{equation}
( - )\bbf\c \cTop \to \cbfTop,   \qq   X  \mapsto  X\bbf,
    \label{2.3.2} \end{equation}
gives the {\em generated border flexible space} $ X\bbf$, with the least c-structure containing the border 
restrictions of the original c-paths, as specified in \ref{2.1}(c). 

	Maps stay `unchanged'; the unit is the reshaping $ X \to X\bbf$. Obviously $ X $ and $ X\bbf $ 
generate the same d-space. 

	Again, border flexible c-spaces are closed in $ \cTop $ under limits. All colimits can be obtained taking 
the colimit in $ \cTop $ and applying the reflector. Sums are preserved by the inclusion and pushouts need 
not be.
 
\Ndt (c) Controlled spaces with a total path-support form a full coreflective subcategory of $ \cTop$, 
 closed under colimits and products.
}
\begin{proof}
(a) The reshaping $ X \to X\pf $ is a universal arrow from $ X $ to the inclusion $ \cpfTop \to \cTop$, which 
creates all limits.

	Closure under sums is trivial, while this does not work with pushouts, as we have seen in \ref{2.2}(b). 
All colimits can be constructed as specified above; or directly, by final structures.

\Ndt (b) The c-space $ X $ is indeed finer than $ X\bbf$: if $ a' = a\rho $ is a border restriction of a c-path of 
$ X$, we can reconstruct $ a = a'\si $ by a global reparametrisation of $ a' $ that brings back the delays.

	Technically, $ a $ is constant on two intervals $ [1, t_1] $ and $ [t_2, 1]$, with $ t_1 < t_2$, and we let:  
    \begin{equation}
\rho(t)  =  t_1 + (t_2 - t_1)t,   \q   \si(t)  =  ((t - t_1)/(t_2 - t_1) \jo 0) \me 1,
    \label{2.3.3} \end{equation}
so that $ a\rho\si = a$. The reshaping $ X \to X\bbf $ is a universal arrow from $ X $ to the inclusion 
$\cbfTop \to \cTop$, and the rest works as in (a), taking into account example \ref{2.2}(c) for pushouts.

\Ndt(c) The coreflector takes a c-space $ X $ to its path-support $ |X|_1 $ equipped with the same 
c-paths. Moreover the path-support functor $ | - |_1\c \cTop \to \Top $ preserves products.
\end{proof}
%

\subsection{Lemma}\label{2.4}
{\em
(a) Let $ a $ be a c-path of the d-space $ \hat{X} $ generated by the c-space $X$. We suppose that one 
of these conditions is satisfied:

\LL

(i) $\; a $ is not constant,

(ii) $ a $ is constant at a point $ x $ of the path-support $ |X|_1$.

\LB
\skp

	Then $ a $ is a restriction of some c-path $b\c x_0 \to x_1$ of $\hat{X}$ between flexible points 
of $X$. We can choose $b$ so that $a$ is its {\em middle restriction}, on $ [1/3, 2/3]$.

\Ndt (b) If the c-space $X$ is preflexible, every c-path of $\hat{X}$ satisfying (i) or (ii) is the middle 
restriction of some c-path of $X$. 
}
\begin{proof}
We prove (a), which trivially implies (b). In case (i) a non-constant c-path $a\c x' \to x''$ of the d-space $\hat{X}$ 
is a general concatenation of c-paths $ a_1, ..., a_n $ which are restrictions of c-paths $\, b_1, ..., b_n \,$ of 
$ X$. There is thus some c-path $ b'\c \!$ $x_0 \to x' $ of $ \hat{X} $ starting from a flexible point of $ X$: either 
a restriction of $b_1$, as in the left figure below
%
    
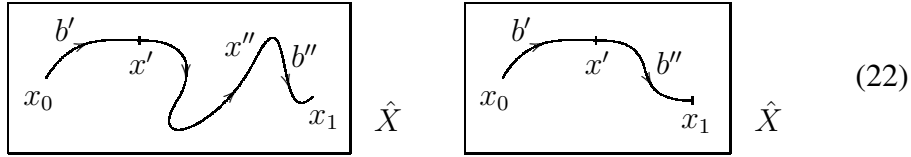
\begin{equation} 
\xy <.5mm, 0mm>:
%
(0,25) *{}; (0,-25) *{};
(90,-10) *{\hat{X}}; (190,-10) *{\hat{X}};
(-2,-6) *{x_0}; (5,13) *{b'}; (25,4) *{x'}; 
(51,9) *{x''}; (68,6) *{b''}; (73,-11) *{x_1}; 
(118,-6) *{x_0}; (125,13) *{b'}; (145,4) *{x'}; (164,4) *{b''}; (171,-12) *{x_1}; 
@i@={(-10,-20), (80,-20), (80,20), (-10,20)},
s0="prev"  @@{;"prev";**@{-}="prev"};
@i@={(110,-20), (180,-20), (180,20), (110,20)},
s0="prev"  @@{;"prev";**@{-}="prev"};
(0,0); (25,10) **\crv{(7,12)&(17,10)}, (25,10); (35,-5) **\crv{(35,10)&(40,5)},
(35,-5); (55,5) **\crv{(25,-20)&(45,-15)}, (55,5); (70,-5) **\crv{(65,25)&(60,-15)},
(120,0); (145,10) **\crv{(127,12)&(137,10)}, (145,10); (170,-6) **\crv{(165,10)&(150,-6)}, 
\POS(8,8) \ar+(3,1.2), \POS(36.7,3) \ar+(0,-3), \POS(48,-5) \ar+(2,2), \POS(63,1) \ar+(1,-3.5),
\POS(128,8) \ar+(3,1.2), \POS(157.7,1) \ar+(1.6,-3),
\POS(-1,0) \arl+(1.5,-1), \POS(24.5,9) \arl+(0,2), 
\POS(34,-5) \arl+(2,-1), \POS(54,5) \arl+(2,-1), \POS(69,-4) \arl+(1,-1.5), 
\POS(119,0) \arl+(1.5,-1), \POS(144.5,9) \arl+(0,2), \POS(170,-7) \arl+(0,2),
\endxy
    \label{2.4.1} \end{equation}
or a trivial loop at $x_0 = x'$, if $x'$ is a flexible point of $ X$. Symmetrically there is in 
$\hat{X}$ a c-path $b''\c x'' \to x_1$ reaching a flexible point of $X$.

	Now, $ a $ is the middle restriction of the regular concatenation $ b = b' * a * b''\c x_0 \to x_1$.
	
	Case (ii) is obvious: there is a c-path $ b\c x_0 \to x_1 $ {\em of} $X$ whose image contains $x$ 
and we can take in $ \hat{X} $ suitable restrictions $\, b'\c x_0 \to x \,$ and $\, b''\c$ $x \to x_1$ of $b$, as in 
the right figure above; or we take two trivial loops, if $x$ is already a flexible point of $X$. Again, $a$ is 
the middle restriction of the regular concatenation $b' * a * b''\c x_0 \to x_1$ in $\hat{X}$.
\end{proof}
%

\subsection{Theorem {\rm (Flexibility and cartesian products)}}\label{2.5}
{\em
(a) The flexible c-space $ (\Pro_i X_i)\hat{\,}$ generated by a product is finer than $ \Pro_i \hat{X}_i$, and 
can be strictly finer.

\Ndt (b) For any family $ (X_i) $ of preflexible c-spaces with a total path-support (see \ref{2.1}(d)):
}
    \begin{equation}
(\Pro X_i)\hat{\,} \; = \; \Pro \hat{X}_i.
    \label{2.5.1} \end{equation}
\begin{proof}
(a) The flexible structure of $ \Pro \hat{X}_i $ contains the structure of $ \Pro X_i$, and also of 
$ (\Pro X_i)\hat{\,}$. The inclusion can be strict, as shown in I.2.7(b) for the products $ X \ti \cI $ or 
$ X \ti \uI$, where $ X = D_c\sing $ has an empty path-support. One could also use the diagonal c-space 
of \ref{2.2}(b), which is not preflexible. 

\Ndt (b) We prove that, if all $ X_i $ are preflexible c-spaces with a total path-support, every c-path 
$a = \lan a_i \ran $ of $ \Pro \hat{X}_i$ is also controlled in $ (\Pro X_i)\hat{\,}$.

	By Lemma \ref{2.4}(b), each $ a_i $ is the middle restriction $ b_i\rho $ of some c-path $ b_i $ of 
$ X_i$, always applying the strictly increasing affine map $ \rho\c \bbI \to \bbI $ with image $ [1/3, 2/3]$. 
Therefore $ a = \lan b_i \ran \rho \,$ is the restriction of a c-path of the product, and belongs to the structure 
$ (\Pro X_i)\hat{\,}$.
\end{proof}
%

\subsection{Corollary}\label{2.6}
{\em
(a) If $ X $ is a preflexible c-space with a total path-support
    \begin{equation}
(X \ti \cI)\hat{\,}  =  \hat{X} \times \uI,   \q\;\;\;   (X \times \uI)\hat{\,}  =  \hat{X} \ti \uI.
    \label{2.6.1} \end{equation}

\ndt (b) In particular: $(\cI \ti \uI)\hat{\,} = \uI^2$.
}

\section{Elementary homotopy theory of d-spaces}\label{s3}

	We recall here the elementary part of homotopy theory in the category $\dTop$ of d-spaces and the 
construction of their fundamental category \cite{G3}, which will be later extended to c-spaces. Some new 
results on c-spaces are already inserted in Proposition \ref{3.3} and Theorem \ref{3.9}.

\subsection{Directed homotopy}\label{3.1}
Homotopy in $\dTop$ is based on the directed interval $ \uI$ and the reversor endofunctor 
$ R\c \dTop \to \dTop$.

	Inside this theory, a map $ a\c \uI \to X $ is simply called a {\em path} -- or a directed path when we 
want to stress the difference with the paths of the underlying space $ UX$. Homotopies are represented 
by maps $ \ph\c X \ti \uI \to Y$, defined on the (directed) cylinder $ X \ti \uI$. This works by a complex 
structure on the interval and the cylinder functor (developed in \cite{G3}, Chapters 1 and 4), of which we 
recall here the initial part. 

	The first-order structure of the interval in $\dTop$ consists of four maps: two {\em faces} 
$ (\ddm, \ddp)$, a {\em degeneracy} $ (e) $ and a {\em reflection} $(r)$
    \begin{equation} \begin{array}{cccc}
\dd^\al\c \sing  \rlr  \uI \cc e,  &\;\;\;\;&  r\c \uI \to \uI\op  & \q (\al = 0, 1),
\\[5pt]
\dd^\al(*)  =  \al,  \;   e(t)  =  *,   &&   r(t)  =  1 - t.
    \label{3.1.1} \end{array} \end{equation}
\begin{small}

	(The same structure works in $ \Top$, using the euclidean interval $ \bbI $ and a trivial reversor $R$, 
the identity of the category.)

\end{small}

\skp	The {\em cylinder} endofunctor $ I_d $ (written as $ I $ if it is clear that we are working in $\dTop$)
    \begin{equation}
I_d\c \dTop \to \dTop,   \qq   I_d  =  - \ti \uI,
    \label{3.1.2} \end{equation}
inherits from this structure four natural transformations, with the same names and symbols:
    \begin{equation} \begin{array}{c}
\dd^\al\c 1  \rlr  I \cc e,   \q\;\;\;   r\c IR \to RI    \q\;\;\;   (\al = 0, 1),
\\[5pt]
\dd^\al(x)  =  (x, \al),   \q   e(x, t)  =  x,   \q   r(x, t)  =  (x, 1 - t).
    \label{3.1.3} \end{array} \end{equation}
\begin{small}

	The component $\, \dd^\al X = X \ti \dd^\al $ on the d-space $ X $ is simply written as $ \dd^\al$, when 
this is not ambiguous; similarly for the other natural transformations. 

\end{small}

	These natural transformations satisfy the identities
    \begin{equation} \begin{array}{ccc}
e\dd^\al  =  1\c  \id \to \id,   &\;\;&   (RrR)r  =  1\c  IR \to IR,
\\[5pt]
(Re)r  =  eR\c  IR \to R,   &&   r(\ddm R)  =  R\ddp\c  R \to RI.
    \label{3.1.4} \end{array} \end{equation}

	Since $ RR = 1$, $ r $ is invertible with $ r^{-1} = RrR\c RI \to IR$. Moreover $r(\ddp R) = R\ddm$.

\vskip 3pt

	A (directed) {\em homotopy} $ \ph\c f \to g\c X \to Y $ of d-spaces is defined as a d-map 
$ \ph\c IX \to Y $ with faces $ \ph\ddm = f $ and $ \ph\ddp = g$. When we want to distinguish the homotopy 
from the {\em map} $ IX \to Y $ which represents it, we write the latter as $\hat{\ph}$. A path in $ X $ is the 
same as a homotopy $ a\c x \to y\c \sing \to X $ between its endpoints, identified to maps $ \sing \to X$.

	Each map $ f\c X \to Y $ has a {\em trivial} endohomotopy, $ e_f\c f \to f$, represented by 
$ f(eX) = (eY)If\c IX \to Y$.

\vskip 2pt

	Every homotopy $ \ph\c f \to g\c X \to Y $ has a {\em reflected homotopy} between the opposite 
d-spaces
    \begin{equation} \begin{array}{c}
\ph\op\c g\op \to f\op\c X\op \to Y\op,
\\[5pt]
(\ph\op)\hat{\,}  \, =  \, R(\hat{\ph})r\c IRX \to RIX \to RY,
    \label{3.1.5} \end{array} \end{equation}
and $\, (\ph\op)\op = \ph$, $\; (e_f)\op = e_{f\op}$.

\vskip 2pt

	There is a {\em whisker composition} for maps and homotopies
    \begin{equation} \begin{array}{c} 
    \xymatrix  @C=30pt @R=20pt
{
~X'~   \ar[r]^-h    &   ~X~   \ar@<7pt>[r]^-f  \ar@<-4pt>[r]_-g  \arv@<2pt>[r]|-{\dw\sst{\ph}} &   
~Y~   \ar[r]^-k    &   ~Y'~   
}
    \label{3.1.6} \end{array} \end{equation}
\vskip-5pt
    \begin{equation*} \begin{array}{c}
k \ci \ph \ci h\c  kfh \to kgh\c  X' \to Y',
\\[3pt]
(k \ci \ph \ci h)\hat{\,}  =  (k \hat{\ph}) (Ih)\c  IX' \to Y',
    \label{3.1.6bis} \end{array} \end{equation*}
\vskip 2pt
\ndt  which will also be written as $ k\ph h$. This ternary operation satisfies obvious axioms of associativity 
and identities (cf.\ \cite{G3}, 1.2.3).


\subsection{Concatenating paths and homotopies}\label{3.2}
In $\dTop$ the {\em standard concatenation pushout} -- pasting two copies of the d-interval, one after the 
other -- can be realised as $ \uI $ itself, with embeddings $c^-, c^+$ covering the first or second half of the 
interval
%
    \begin{equation} \begin{array}{c} 
    \xymatrix  @C=11pt @R=2pt
{
~\sing~   \ar[rr]^-{\ddp}   \ar[dd]_(.45){\ddm}    &&   ~\uI~    \ar[dd]^(.45){c^-}   &&  c^-(t)  =  t/2,
\\ 
&&   \ar@{.}@/_/[ld]          &&   \q\;   c^+(t)  =  (t + 1)/2.
\\ 
~\uI~   \ar[rr]_-{c^+}    &&   ~\uI~
}
    \label{3.2.1} \end{array} \end{equation}

	Indeed, given two consecutive paths $ a, b\c \uI \to X $ (with $ a\ddp= b\ddm$), their concatenation 
$ a * b $ is a map $ \uI \to X$.

	More generally we want to concatenate two consecutive homotopies $ \ph\c$ $f \to g\c X \to Y $ 
and $ \psi\c g \to h\c X \to Y$, with $ \ph(\ddp X) = g = \psi(\ddm Y)$, defining
    \begin{equation} \begin{array}{c}
\chi  =  \ph * \psi\c f \to h\c X \to Y,
\\[5pt]
\chi(x, t)  =  \ph(x, 2t),   \q\,   \for  0 \le t \le 1/2,
\\[3pt]
\chi(x, t)  =  \psi(x, 2t - 1), \;\;  \for  1/2 \le t \le 1.
    \label{3.2.2} \end{array} \end{equation}

	As in $ \Top$, we do obtain a map $ \chi\c X \ti \uI \to Y$, because of the following proposition.

\subsection{Proposition {\rm(Concatenating homotopies)}}\label{3.3}
{\em
(a) For every d-space $ X$, the functor $ X \ti -\c \dTop \to \dTop $ preserves the standard concatenation 
pushout \eqref{3.2.1}, giving the {\em concatenation pushout of the cylinder functor} $ I = - \ti \uI $
%
    \begin{equation} \begin{array}{c} 
    \xymatrix  @C=8pt @R=3pt
{
~X~   \ar[rr]^-\ddp   \ar[dd]_(.45)\ddm    &&   ~X \ti \uI~    \ar[dd]^(.45){c^-}   &  c^-(x, t)  =  (x, t/2),
\\ 
&&   \ar@{.}@/_/[ld]          &   \q\;   c^+(x, t)  =  (x, (t+1)/2).
\\ 
~X \ti \uI~   \ar[rr]_-{c^+}    &&   ~X \ti \uI~
}
    \label{3.3.1} \end{array} \end{equation}

\Ndt (b) If $ X $ is a {\em flexible} c-space, this is also a pushout in $ \cTop$, so that the functor 
$ X \ti -\c \cTop \to \cTop $ preserves the pushout \eqref{3.2.1}.
}
\begin{proof}
(a) We copy the proof of \cite{G3}, 1.4.9, in a more detailed way that will be used in the next section.

	The square \eqref{3.3.1} becomes a pushout in $\Top$, because $UX \ti [0, 1/2]$ and $UX \ti [1/2, 1]$ 
form a finite closed cover of $ UX \ti \bbI$, so that each mapping defined on the latter and continuous on 
these closed parts is continuous.

	Consider then the map $ \chi\c UX \ti \bbI \to UY $ obtained by pasting two d-maps 
$ \ph, \psi\c X \ti \uI \to Y $ on the topological pushout $ UX \ti \bbI$, as in \eqref{3.2.2}. We want to prove 
that $ \chi $ is a d-map $ X \ti \uI \to Y$.

	Let $\lan a, h \ran\c \uI \to IX $ be any d-map, with $ a\c \uI \to X $ and $  h\c \uI \to \uI$.

\Ndt(i) If the image of $ h $ is contained in the first (resp.\ second half) of $ [0, 1]$, then $ \chi \lan a, h \ran$ is 
directed, because
    \begin{equation} \begin{array}{c}
\chi(a(t), h(t))  =  \ph(a(t), 2h(t))
\\[5pt]
(\mbox{resp.} \;   \chi(a(t), h(t))  =  \ph(a(t), 2h(t) - 1)),  \,
    \label{3.3.2} \end{array} \end{equation}
and, in this case, $ 2h $ (resp.\ $ 2h - 1$) is a map $ \uI \to \uI $ (just an increasing continuous mapping).

\Ndt (ii) Otherwise, there is some $ t_1 \! \in \; ]0, 1[ $ such that $ h(t_1) = 1/2$, and we can assume that 
$ t_1 = 1/2 $ (up to pre-composing with an automorphism of $ \uI$). Now, the path 
$ \chi\lan a, h \ran\c \bbI \to UY $ is directed in $ Y$, because it is the concatenation of the following 
two directed paths $ c_i\c \uI \to Y $
    \begin{equation} \begin{array}{c}
c_1(t)  =  \chi(a(t/2), h(t/2))  =  \ph(a(t/2),  2h(t/2)),
\\[5pt]
c_2(t)  =  \chi(a((t + 1)/2), h((t+1)/2))
\\[3pt]
\qq\;\;  =  \psi(a((t + 1)/2),  2h((t + 1)/2) - 1).
    \label{3.3.3} \end{array} \end{equation}

	Note that {\em we are using the splitting property of} $ a$, in the d-space $X$.

\Ndt (b) A straightforward consequence: by hypothesis $X$ is a d-space, and the pushout \eqref{3.3.1} is 
preserved by the embedding in $ \cTop$.
\end{proof}
%

\subsection{Homotopies of paths}\label{3.4}
Operating on the standard concatenation pushout \eqref{3.2.1} with the functors $ - \ti \uI $ and $ \uI \ti -$, 
we get the following pushouts
%
    \begin{equation} \begin{array}{c} 
    \xymatrix  @C=8pt @R=10pt
{
~\uI~   \ar[rr]^-{\ddp \ti 1}   \ar[dd]_(.4){\ddm \ti 1}    &&   ~\uI \ti \uI~    \ar[dd]^(.4){\Phi'}  &&&
~\uI~   \ar[rr]^-{1 \ti \ddp}   \ar[dd]_(.4){1 \ti \ddm}    &&   ~\uI \ti \uI~    \ar[dd]^(.4){\Phi''}
\\ 
&&   \ar@{.}@/_/[ld]    &&&&&   \ar@{.}@/_/[ld]         
\\ 
~\uI \ti \uI~   \ar[rr]_-{\Psi'}    &&   ~\uI \ti \uI~   &&&
~\uI \ti \uI~   \ar[rr]_-{\Psi''}    &&   ~\uI \ti \uI~
}
    \label{3.4.1} \end{array} \end{equation}
that we use now to concatenate double paths, horizontally and vertically.

	A {\em double path} in the d-space $ X $ is a d-map $ H\c \uI^2 \to X$. Its four faces are paths 
in $ X$, between four points, its vertices
%
    \begin{equation} 
\xy <.5mm, 0mm>:
%
(0,25) *{}; (140,-25) *{};
(25,0) *{H}; (0,0) *{\sst{H\dd^-_1}}; (50,0) *{\sst{H\dd^+_1}}; 
(25,-22) *{\sst{H\dd^-_2}}; (25,22) *{\sst{H\dd^+_2}}; 
(125,-15) *{\sst{s}}; (119,0) *{\sst{t}};
@i@={(10,-15), (40,-15), (40,15), (10,15)},
s0="prev"  @@{;"prev";**@{-}="prev"};
\POS(110,-10) \ar+(20,0),  \POS(115,-15) \ar+(0,20),
\endxy
    \label{3.4.2} \end{equation}
    \begin{equation*} \begin{array}{ccc}
\dd^\al_1  =  \dd^\al \ti \uI\c \uI \to \uI^2,   &\q&   \dd^\al_1(t)  =  (\al, t),
\\[5pt]
\dd^\al_2  =  \uI \ti \dd^\al\c \uI \to \uI^2,   &&   \dd^\al_2(s)  =  (s, \al).
    \label{3.4.2bis} \end{array} \end{equation*}

	Given two horizontally consecutive double paths $ H, K\c \uI \ti \uI \to X $ ($H\dd^+_1 = K\dd^-_1$), 
the left pushout \eqref{3.4.1} gives a map $ H *_1 K\c \uI \ti \uI \to X$, the horizontal concatenation of 
$ H $ and $ K$.

	Symmetrically, given two double paths $H, L\c \uI \ti \uI \to X$ which are vertically consecutive, the right 
pushout \eqref{3.4.1} gives the vertical concatenation $ H *_2 L\c \uI \ti \uI \to X$.

\begin{small}
\vskip 2pt

	It is evident and well-known that these operations satisfy -- even strictly -- the middle-four interchange 
property:
\vskip-8pt

$$ (H *_1 K) *_2 (L *_1 M) \, = \, (H *_2 L) *_1 (K *_2 M).$$
\vskip-2pt

In $ \Top $ this induces the commutativity of $\pi_2$ (using double paths with degenerate faces).

\end{small}
\vskip 2pt

	A {\em path-homotopy with fixed endpoints}, denoted as
    \begin{equation}
H\c a \to_2 b\c x \to y,
    \label{3.4.3} \end{equation}
is a double path $ H\c \uI^2 \to X $ whose vertical faces $H\dd^\al_1$ are trivial loops at two points 
$ x, y $ (constant on the thick edges) 
%
    \begin{equation} 
\xy <.4mm, 0mm>:
%
(0,25) *{}; (190,-25) *{};
(25,0) *{H}; (2,0) *{e_x}; (48,0) *{e_y}; (25,-21) *{a}; (25,22) *{b}; 
(130,10) *{H(0, t) = x,   \;\;    H(1, t) = y,}; (130,-10) *{H(-, 0)  =  a,   \;\;   H(-, 1)  =  b.};
@i@={(10,-15), (40,-15), (40,15), (10,15)},
s0="prev"  @@{;"prev";**@{-}="prev"};
 \POS(10.4,-15) \arl+(0,30),  \POS(39.6,-15) \arl+(0,30),
\endxy
    \label{3.4.4} \end{equation}

	$H $ will also be called a {\em 2-path}, by analogy with the 2-cells of a 2-category. The existence of
 $H\c a \to_2 b $ gives a preorder relation $ a \precc_2  b $ consistent with concatenation, because 2-paths 
 are obviously closed under the vertical and horizontal concatenation of double paths.
 
	We now use the {\em 2-equivalence} relation $a' \eqq a''$ spanned by the preorder $\precc_2$: 
 there exists a finite sequence of paths $a' \precc_2 a_1  \succeq_2 a_2 \, ...  \precc_2 a''$ (between the 
 same endpoints). It is also consistent with concatenation.
 
	A class of paths $ [a] $ {\em up to 2-equivalence} is a class of this equivalence relation. These 
 classes become now the arrows of the fundamental category.

\subsection{Definition {\rm  \cite{G1, G3}}}\label{3.5}
In the {\em fundamental category} $\uPi_1(X)$ of a d-space $X$: 

\LL

\ndt - the vertices are the points of $ X$,

\Ndt - the arrows are the equivalence classes $ [a]\c x \to y $ of paths $ a\c x \to y$, up to the 2-equivalence 
relation $ a' \eq\!_2 \, a'' $ described above (spanned by homotopies with fixed endpoints $ a' \to_2 a''$),

\Ndt - the composition is induced by concatenation of consecutive paths, and the identities are induced by 
trivial loops
    \begin{equation}
[a][b]  =  [a * b],   \q\;\;\;   1_x  =  [e_x].
    \label{3.5.1} \end{equation}

\LB

	The operation is well defined, because $\, \eq\!_2 $ is consistent with path concatenation. 
Associativity and identities follow from Lemma \ref{3.7}, which also implies that a global (increasing) 
reparametrisation $ \rho\c \bbI \to \bbI \, $ gives $ [a] = [a\rho]$.

	A d-map $ f\c X \to Y $ produces a functor
    \begin{equation} \begin{array}{c}
f_*  =  \uPi_1(f)\c \uPi_1(X) \to \uPi_1(Y),
\\[5pt]
f_*(x)  =  f(x),   \q   f_*[a]  =  [fa].
    \label{3.5.2} \end{array} \end{equation}

\begin{small}

	In fact a 2-path $H\c a' \to_2 a''$ gives a 2-path $fH\c fa' \to_2 fa''$, and a concatenation 
$ a' * a'' $ of paths in $X$ is sent to a concatenation $fa' * fa''$ in $ Y$.

\end{small}

\skp
	We have thus a functor with values in the category of small categories
    \begin{equation}
\uPi_1\c \dTop \to \Cat.
    \label{3.5.3} \end{equation}

\begin{small}

	The next theorem refers to the directed homotopy structure of $\Cat$, based on the 
directed interval $\two$  (\cite{G3}, 1.2.2, 4.3.2); in this light, a natural transformation is a 
directed homotopy of functors.

\end{small}

\subsection{Theorem {\rm (Homotopy invariance, I)}}\label{3.6}
{\em
(a) The functor $ \uPi_1 $ is homotopy invariant, in the sense that a (directed) homotopy $ \ph\c f \to g\c X \to Y $ induces a natural transformation
    \begin{equation} \begin{array}{c}
\ph_*\c f_* \to g_*\c \uPi_1(X) \to \uPi_1(Y),
\\[5pt]
\ph_*(x)  =  [\ph_x]\c f(x) \to g(x),
    \label{3.6.1} \end{array} \end{equation}
where $\ph_x = \hat{\ph}(x \ti \ul{i})\c \sing \ti \uI \to Y$ is a directed path in $Y$.

\Ndt (b) Induction preserves whisker composition: a whisker composition $\psi = k \ci \ph \ci h\c kfh \to kgh $ 
of homotopies and maps (as in \eqref{3.1.6}) is taken to the corresponding whisker composition of natural 
transformations and functors:
}
    \begin{equation}
\psi_*  =  k_* \ci \ph_* \ci h_*\c  k_*f_*h_* \to k_*g_*h_*\c  \uPi_1(X') \to \uPi_1(Y').
    \label{3.6.2} \end{equation}
\begin{proof}
(a) Again we write out the proof of \cite{G3}, Theorem 3.2.4(c), to extend it later to the present 
setting.

	Let $ a\c x \to x' $ be a path in $ X$, represented by a map $ a\c \uI \to X$. The homotopy 
$ \ph\c X \ti \uI \to Y $ gives a double path $ \ph_a $ in $ Y$, with the following faces
%
    \begin{equation} \begin{array}{c} 
    \xymatrix  @C=7pt @R=3pt
{
~g(x)~   \ar[rr]^-{ga}    &&   ~g(x')~   
\\ 
&  \ph_a  &&  \q  \ph_a  = \hat{\ph}(a \ti \ul{i})\c \uI \ti \uI \to Y,
\\ 
~f(x)~   \ar[rr]_-{fa}   \ar[uu]^(.45){\ph_x}    &&   ~f(x')~  \ar[uu]_(.45){\ph_{x'}}
}
    \label{3.6.3} \end{array} \end{equation}
where the vertical direction agrees with the present graphic representation of double paths.
We have to prove that this double path induces a commutative square in $ \uPi_1(Y)$: 
$f_*(a) * \ph_*(x')  =  \ph_*(x) * g_*(a).$

	This comes out nearly for free from the higher structure of the directed interval $ \uI $ (as a 
{\em dioid} \cite{G3}, 1.1.7), which includes the $\max$ and $\min$ operations $g^\al$ (with unit $\dd^\al$)
    \begin{equation}
g^\al\c \uI^2 \to \uI,   \q\;   g^-(s, t)  =  s \jo t,  \;\;\;\;  g^+(s, t)  =  s \me t,	
    \label{3.6.4} \end{equation}
%
%
    \begin{equation*} 
\xy <.5mm, 0mm>:
%
(0,27) *{}; (170,-25) *{};
(25,0) *{g^-}; (5,0) *{\sst{\ul{i}}}; (47,0) *{\sst{e_1}}; (25,-22) *{\sst{\ul{i}}}; (25,21) *{\sst{e_1}}; 
(125,0) *{g^+}; (104,0) *{\sst{e_0}}; (145,0) *{\sst{\ul{i}}}; (125,-22) *{\sst{e_0}}; (125,21) *{\sst{\ul{i}}}; 
@i@={(10,-15), (40,-15), (40,15), (10,15)},
s0="prev"  @@{;"prev";**@{-}="prev"};
@i@={(110,-15), (140,-15), (140,15), (110,15)},
s0="prev"  @@{;"prev";**@{-}="prev"};
\POS(10,14.5) \arl+(30,0),  \POS(39.6,-15) \arl+(0,30),
\POS(110,-14.5) \arl+(30,0),  \POS(110.4,-15) \arl+(0,30),
\endxy
    \label{3.6.4bis} \end{equation*}

	Working with the double path \eqref{3.6.3} and the new ones, by horizontal concatenation, we get a 
double path $\, (\ph_x g^+ *_1 \ph_a) *_1 \ph_{x'}g^- \,$ (the parentheses are irrelevant)
%
    \begin{equation} \begin{array}{c} 
    \xymatrix  @C=6pt @R=3pt
{
~f(x)~   \ar[rr]^-{\ph_x}    &&   ~g(x)~   \ar[rr]^-{ga}    &&   ~g(x')~   \ar[rr]^-{e_{gx'}}    &&   ~g(x')~   &    ~\,~   
\\ 
&  \ph_x g^+   &&   ~~\ph_a~~    &&   \ph_{x'} g^-
\\ 
~f(x)~   \ar[rr]_-{e_{fx}}   \ar[uu]^(.45){e_{fx}}    &&   ~f(x)~  \ar[rr]_-{fa}  \ar[uu]_(.45){\ph_{x}}  &&
~f(x')~  \ar[rr]_-{\ph_{x'}}  \ar[uu]^(.45){\ph_{x'}}  &&  ~g(x')~   \ar[uu]_(.45){e_{gx'}}
}
    \label{3.6.5} \end{array} \end{equation}

	It is a 2-path from $\, (e_{fx} * fa) * \ph_{x'} \, $ to $\, (\ph_x * ga) * e_{gx'}\, $ in $Y$, and proves our 
identity in $\uPi_1(Y)$
$$
f_*(a) * \ph_*(x') =  [fa][\ph_{x'}]  =  [\ph_x][ga]  =  \ph_*(x) * g_*(a).
$$

\Ndt (b) The component $ \psi_*(x) = [\psi_x] $ at $ x \in X' $ is the 2-equivalence class of the following 
d-path $kfh(x) \to kgh(x)$ in $Y'$:
$$
\psi_x  = \hat{\psi}(x \ti \ul{i}) = k\hat{\ph} (Ih)(x \ti \ul{i}) = k\hat{\ph}(h(x) \ti \ul{i}) = k \ph_{hx}.
$$
and coincides with the component $(k_* \ci \ph_* \ci h_*)(x)$.
\end{proof}
%

\subsection{Lemma {\rm (Reparametrisation and homotopy)}}\label{3.7}
{\em
Let $ a\c x \to y $ be a path in the d-space $ X$.

\Ndt (a) If $ \rho, \si\c \bbI \to \bbI $ are global reparametrisations and $ \rho \le \si$, there exists a 
2-path 
    \begin{equation}
H\c a\rho \to_2 a\si\c \uI \to X,    \q\;\;    H\c \uI^2 \to X.   \q\;\;
    \label{3.7.1} \end{equation}

\Ndt (b) In particular, using the global reparametrisations \eqref{1.5.1} and  \eqref{1.5.2}, 
we get 2-paths of {\em acceleration} and {\em associativity}: 
    \begin{equation}
e_x * a \to_2 a,   \qq\;\;   a \to_2 a * e_y,
    \label{3.7.2} \end{equation}
    \begin{equation}
a * (b * c) \to_2 a * b *c,   \q   a * b * c \to_2 (a * b) * c.
    \label{3.7.3} \end{equation}

\Ndt (d) For any global reparametrisation $\, \rho\c \bbI \to \bbI \,$ there are global reparametrisations 
$ \rho', \rho'' $ and 2-paths
    \begin{equation}
a\rho' \to_2 a \to_2 a\rho'',   \q   a\rho' \to_2 a\rho \to_2 a\rho'',
    \label{3.7.4} \end{equation}
proving that $a \, \eqq \, a\rho$.
}
\begin{proof}
(a) We build a 2-path $\, \Phi\c \rho \to_2 \si\c \uI \to \uI$
    \begin{equation}
\Phi\c \uI^2  \to \uI,   \q  \Phi(s, t)  =  (1 - t)\rho(s) + t\si(s),
    \label{3.7.5} \end{equation}
as the affine interpolation from $ \rho $ to $ \si$, which is increasing because $ \rho \le \si$. 
(This corresponds to a particular case, where $ a $ is the identity of $ \uI$.)

	Then, composing $ \Phi $ with the path $ a\c \uI \to X$, we get a 2-path
    \begin{equation}
a\Phi\c a\rho \to_2 a\si\c \uI \to X.
    \label{3.7.6} \end{equation}

\ndt (b) Follows from (a), using the global reparametrisations $ \rho \le \ul{i} \le \si$ of \eqref{1.5.1} 
and \eqref{1.5.2}.

\Ndt (c) Follows from (a), using $ \rho' = \rho \me \ul{i} \,$ and $ \rho'' = \rho \jo \ul{i}$.
\end{proof}
%

\subsection{The cocylinder functor}\label{3.8}
Coming back to first-order homotopy theory in $\dTop$, we recall that the cylinder functor 
$ I = - \ti \uI\c \dTop \to \dTop $ has a right adjoint, the {\em cocylinder functor}, or {\em path functor}
    \begin{equation}
P\c \dTop \to \dTop,   \q   P(Y)  =  Y^\uI   \q\;\;   (I \adj P).
    \label{3.8.1} \end{equation}

	In fact, the d-interval $ \uI $ is exponentiable in $\dTop$, by Theorem \ref{3.9} (see below). The 
d-space $P(Y)$ is the set of d-paths $ \dTop(\uI, Y) $ equipped with the compact-open topology (induced 
by the topological path-space $\, P|Y| = \Top(\bbI, |Y|)) $ and the d-structure described in the statement 
of \ref{3.9}.

	A homotopy $ \ph\c f^- \to f^+\c X \to Y $ can thus be equivalently defined

\LL

\ndt - by a map $\, \hat{\ph}\c X \ti \uI \to Y \, $ with $\, \hat{\ph}\dd^\al = f^\al$, for $ \al = \pm$,

\ndt - by a map $\check{\ph}\c X \to Y^\uI $ with $ \dd^\al\check{\ph} = f^\al $ (with respect to the faces 
$ \dd^\al\c P \to 1 $ specified below).
	
\LB
\vskip 2pt

	The path functor is equipped with a first-order structure dual to that of the cylinder functor, with the 
same reversor $ R$; it consists of four natural transformations, still called {\em faces}, 
{\em degeneracies} and {\em reflection}, and denoted by the same symbols (for $a\c \uI \to X$)
    \begin{equation} \begin{array}{ccc}
\dd^\al\c P \lrw 1 \cc e,   &&   \dd^\al(a)  =  a\dd^\al,   \q   e(a)  =  ae,
\\[5pt]
r\c RP \to PR,   &\;\;\;&   r(a)  =  (Ra)r\c \uI  \to \uI\op \to X\op.
    \label{3.8.2} \end{array} \end{equation}

	In $ \cTop $ this part will have a complex reworking: the interval $ \cI $ is not exponentiable, as
will be shown in \ref{4.7}(b). On the other hand, the interval $\uI$ is still exponentiable in $\cTop$, as 
we prove in \ref{3.9}(b), and gives rise to a {\em functor of flexible paths}, that represents flexible 
homotopies (in \eqref{4.1.4}).

\subsection{Theorem {\rm (Exponentiable spaces)}}\label{3.9}
{\em
Let $ K $ be a d-space on a locally compact support $|K|$.

\Ndt (a) The d-space $ K $ is exponentiable in $\dTop$. For every d-space $ Y$
    \begin{equation}
Y^K  =  \dTop(K, Y)  \sub  \Top(|K|, |Y|)
    \label{3.9.1} \end{equation}
is the set of d-maps $K \to Y$, equipped with:

\LL

\ndt - the compact-open topology (restricted from the topological exponential $|Y|^{|K|})$

\ndt - the d-structure where a path $ c\c \bbI \to |Y^K| \sub |Y|^{|K|}$ is directed if and 
only if the corresponding topological map $\hat{c}\c \bbI \ti |K| \to |Y| $ is a d-map $ \uI \ti K \to Y$.

\LB

\Ndt (b) The d-space $ K $ is also exponentiable in $ \cTop$, extending the definition of the 
previous internal hom: for every c-space $ Y$
    \begin{equation}
Y^K  =  \cTop(K, Y)  =  \dTop(K, \Fl Y)  \sub  \Top(|K|, |Y|),
    \label{3.9.2} \end{equation}
is the set of c-maps $ K \to Y$, with the compact-open topology and the c-structure where a path 
$\, c\c \bbI \to |Y^K| \sub |Y|^{|K|} $ is controlled if and only if the corresponding map 
$\, \hat{c}\c \bbI \ti |K| \to |Y| $ is a c-map $\, \cI \ti K \to Y$. Moreover, all the elements of $\, Y^K $ are 
flexible.

	If $ Y $ is a d-space, we find the same exponential $Y^K$.
}
\begin{proof}
(a) This is Theorem 1.4.8 of \cite{G3}. We rewrite the proof in a more detailed way, to be used in (b).

	In the domain of topological spaces, it is well known that a locally compact space $K$ is 
exponentiable: the space $Y^K$ is the set of maps $\Top(K, Y)$ with the compact-open topology, 
and the adjunction consists of the natural bijection
    \begin{equation}
\hat{\,}\c \Top(X, Y^K) \to \Top(X \ti K, Y),    \q   \hat{f}(x, k)  =  f(x)(k).
    \label{3.9.3} \end{equation}

	Coming back to d-spaces, we verify the axioms (dsp.0)-(dsp.2) for the d-structure of $Y^K$ 
defined above.

\Ndt (i) ({\em Constant paths}) If $c\c \bbI \to |Y|^{|K|}$ is constant at the d-map $g\c K \to Y$, then 
$\, \hat{c} \,$ can be factorised as $gp_2\c \uI \ti K \to K \to Y$, and is a d-map.

\Ndt (ii) {\em Concatenation}) Let $c = c_1 * c_2\c \bbI \to |Y|^{|K|} $, with $\hat{c}_i\c \uI \ti K \to Y$. 
By Proposition \ref{3.3}(a), the product $- \ti K$ preserves the concatenation pushout \eqref{3.2.1}. 
Therefore $ \hat{c}$, as the pasting of $\hat{c}_1, \hat{c}_2$ on this pushout, is a directed map.

\Ndt (iii) ({\em Partial reparametrisation}) For $ \hat{c}\c \uI \ti K \to Y $ and $ h\c \uI \to \uI$, the map 
$ (ch)\hat{\,} = \hat{c}(h \ti K) $ is directed.

\skp	We prove now that \eqref{3.9.3} restricts to a bijection between $ \dTop(X, Y^K) $ 
and $ \dTop(X \ti K, Y)$. 

In fact, we have a chain of equivalent conditions

\LL

\ndt - $f\c X \to Y^K$  is a d-map,

\Ndt - $\, \all a \in X\shp, \;  fa\c \uI \to Y^K$ is a d-path,

\Ndt - $\, \all a \in X\shp, \;  (fa)\hat{\,} = \hat{f}(a \ti K)\c \uI \ti K \to X \ti K \to Y$ is a d-map,

\Ndt - $\, \all a \in X\shp, \, \all h \in \uI\shp, \, \all b \in K\shp,  \;  \hat{f}\lan ah, b\ran \c \uI \to Y$ is a d-map,

\Ndt - $\hat{f}\c X \ti K \to Y$ is a d-map,

\LB

\Ndt taking into account, for the last equivalence, that $\, ah\c \uI \to X $ is an arbitrary d-path.

\Ndt (b) To verify that the c-structure of $ Y^K $ is well formed, point (i) works as above: 
if $ c\c \bbI \to Y^K $ is constant at any c-map $ g\c K \to Y$, the associated map $\, \hat{c} \,$ 
can be factorised as $ gp_2\c \cI \ti K \to K \to Y$, and is a c-map. This also proves that all the 
`points' of $Y^K$ are flexible.

	Point (ii) is also proved as above: we can apply Proposition \ref{3.3}(b), since $ K $ is a 
flexible c-space.

	Point (iii) is now about {\em global reparametrisation}. If $ h\c \uI \to \uI $ is a surjective 
increasing map, $ h \ti K $ is also surjective and $ (ch)\hat{\,} = \hat{c}(h \ti K) $ is controlled.

\vskip 2pt

	To prove that \eqref{3.9.3} restricts to a bijection between $ \cTop(X, Y^K) $ 
and $ \cTop(X \ti K, Y)$ it is sufficient to rewrite the previous chain of equivalences, replacing 
the interval $\uI$ by $\cI$ and the prefix `d' by `c'.

\vskip 2pt

	Finally, if $ Y $ is a d-space, the condition that $\, \hat{c} \,$ be a c-map $ \cI \ti K \to Y $ is 
equivalent to a d-map $ \uI \ti K \to Y$, because $ (\cI \ti K)\hat{\,} = \uI \ti K$, by Theorem \ref{2.5}(b).
\end{proof}
%

\section{Paths and double paths in controlled spaces}\label{s4}

	We now begin to develop homotopy theory in the category $ \cTop$. We work in a concrete, 
naive way, to prepare the definition of the fundamental category of a c-space (in the next section).
This is achieved by a hybrid theory based on general paths, parametrised on $ \cI$, and their 
flexible deformations parametrised on $ \uI$: in other words, we work with the hybrid square 
$ \cI \ti \uI$. 

	Loosely speaking, as made precise in \ref{4.5}(e), if we work on the standard c-square $\cI^2$ our 
homotopies of paths cannot be concatenated; on the other hand, working with the flexible square $\uI^2$ 
we would only get the fundamental category of the flexible part of $X$. But we shall see in Part III that, 
for a border flexible c-space, one can equivalently use the standard square $\cI^2$.

	$X$ and $ Y $ are c-spaces. A controlled path in a c-space is simply called a {\em path}, or a 
{\em general path} if we want to stress that it is not assumed to be flexible; a path in the underlying 
topological space is called a topological path.

\subsection{Paths and homotopies}\label{4.1}
(a) In $ \cTop $ we want to use both the standard c-interval $ \cI$ and the flexible interval $ \uI$. 
They have a similar first-order structure (two faces, a degeneracy and a reflection), already examined 
for $ \uI $ in \eqref{3.1.1}
    \begin{equation} \begin{array}{cccc}
\dd^\al\c \sing  \rlr  \cI \cc e,  &\;\;\;\;&  r\c \cI \to \cI\op,
\\[5pt]
\dd^\al\c \sing  \rlr  \uI \cc e,  &\;\;\;\;&  r\c \uI \to \uI\op.
    \label{4.1.1} \end{array} \end{equation}

	A map $ x\c \sing \to X $ is a flexible point in the c-space $X$; a map $ a\c \cI \to X $ is a c-path 
in $ X$, also called a path, or a general path. 

A map $ a\c \uI \to X $ is a flexible path. The latter is 
{\em flexibly reversible} if also the reversed topological path $ a\shp = ar $ is a flexible path of the 
c-space $ X$; this is equivalent to a c-map $ \bbI^\sim \to X$, defined on the reversible d-interval 
(whose d-structure is generated by $ \ul{i} $ and $ r$, see I.2.4(c)).

\Ndt  (b) A (general) {\em homotopy} $ \ph\c f \to g\c X \to Y $ is a map $ \ph\c X \ti \cI \to Y$, with faces 
$ f = \ph(1 \ti \ddm) $ and $ g = \ph(1 \ti \ddp)$. In particular, a homotopy $ a\c x \to y\c \sing \to X $ is a 
general path (between flexible points), identifying $ \sing \ti \cI $ with $ \cI$.

	A {\em flexible homotopy} $ \ph\c f \to g $ is a map $ \ph\c X \ti \uI \to Y $ with faces as above. A 
flexible homotopy $ a\c x \to y\c \sing \to X $ is the same as a flexible path in $ X$. A flexible homotopy 
$ \ph\c X \ti \uI \to Y $ is {\em flexibly reversible} if also the reflected topological homotopy 
$ \ph\shp = \ph(X \ti r) $ is a flexible homotopy.

	Given a general (resp.\ flexible, flexibly reversible) homotopy $ \ph\c f \to g, $ every flexible point 
$ x \in X $ gives a general (resp.\ flexible, flexibly reversible) path $ \ph(x, -)\c f(x) \to g(x) $ in $ Y$. 
A mere point gives a topological path. (According to the structure of $ X$, it can give something more: 
see \ref{4.4}(b).)

\Ndt (c) The reversor $ R $ transforms a general homotopy $ \ph\c f \to g\c X \to Y $ into the 
{\em reflected} one
    \begin{equation} \begin{array}{c}
\ph\op\c g\op \to f\op\c X\op \to Y\op,		
\\[5pt]
(\ph\op)\hat{\,}  =  R(\hat{\ph})(X\op \ti r)\c  X\op \ti \cI \to (X \ti \cI)\op \to Y\op.
    \label{4.1.2} \end{array} \end{equation}

	(As in \ref{3.1}, the representative map $ IX \to Y $ is written as $ \hat{\ph} $ when it should
be distinguished from the homotopy $ \ph$.) The flexible case works similarly
    \begin{equation} 
(\ph\op)\hat{\,}  =  R(\hat{\ph})(X\op \ti r)\c  X\op \ti \uI \to (X \ti \uI)\op \to Y\op.
    \label{4.1.3} \end{equation}

\Ndt (d) General and flexible homotopies have a whisker composition with maps, as in \eqref{3.1.6}. 
(We shall see that they cannot be concatenated: this is verified in \ref{4.7}(b), (c), for the special case 
$X = \cI$.)

\vskip 3pt
\begin{small}

	More formally, we are working with the {\em standard cylinder functor} $I_c$ and the 
{\em flexible cylinder functor} $I_F$
    \begin{equation} \begin{array}{ccc}
I_c\c \cTop \to \cTop,   &\q&   I_c  =  - \ti \cI,	
\\[5pt]
I_F\c \cTop \to \cTop,   &&   I_F  =  - \ti \uI,
    \label{4.1.4} \end{array} \end{equation}
supplying the category $ \cTop $ with structures that will be investigated in a sequel.

\end{small}

\Ndt (e) Within flexible c-spaces, both general and flexible homotopies coincide with the homotopies 
of d-spaces. The second point is obvious, the first comes from Corollary \ref{2.6}: if $ X $ and $ Y $ 
are flexible, a general homotopy $ \ph\c X \ti \cI \to Y $ amounts to a map defined on 
$ (X \ti \cI)\hat{\,} = X \ti \uI$.

\vskip 2pt

	We now want to form the fundamental category of a c-space.

\subsection{Concatenation of c-paths}\label{4.2}
(a) The concatenation of general paths is formalised in a slightly more complex way than for d-spaces 
(in \ref{3.2}), because of the failure of the path-splitting property.

	The standard concatenation pushout of c-spaces is realised as $\cJ$, the two-jump structure on 
the euclidean interval $ [0, 1] $ (recalled in \ref{1.3}(c))
%
    \begin{equation} \begin{array}{c} 
    \xymatrix  @C=11pt @R=1pt
{
~\sing~   \ar[rr]^-{\ddp}   \ar[dd]_(.45){\ddm}    &&   ~\cI~    \ar[dd]^(.45){c^-}   &&  c^-(t)  =  t/2,
\\ 
&&   \ar@{.}@/_/[ld]          &&   \q\;   c^+(t)  =  (t + 1)/2.
\\ 
~\cI~   \ar[rr]_-{c^+}    &&   ~\cJ~
}
    \label{4.2.1} \end{array} \end{equation}

	Now, given two maps $a', a''\c \cI \to X$ such that $a'\ddp= a''\ddm$, we get a map $a\c \cJ \to X$ 
such that $ ac^- = a'$, $ac^+ = a''$, which `is' also a path $ \cI \to X$, because $ \cI $ has a finer structure. 

	More formally, we can introduce the {\em concatenation map} $ \ka\c \cI \to \cJ$, a reshaping, 
and define $ a' * a'' = a\ka\c \cI \to X$.

\vskip 5pt
\begin{small}

	This procedure is not infrequent in homotopy theory. For instance, chain complexes of abelian 
groups have a similar behaviour: pasting two copies of the interval (or a cylinder) yields a different 
object, related to the former by a concatenation map (cf.\ \cite{G3}, Section \ref{4.4}, where we mostly 
work with the cocylinder). The same happens in $\Cat$, whose two-jump interval is the ordinal $\three$.

\end{small}

\Ndt (b) The problem here is that the pushout \eqref{4.2.1} need {\em not} be preserved by a functor 
$ X \ti -\c \cTop \to \cTop$.

	Trying to adapt the proof of Proposition \ref{3.3} to the category $ \cTop $ {\em and} this pushout, 
we note that:

\LL

\ndt - case (i) still holds: here it can only occur for a c-map $ h\c \cI \to \cI $ constant at $ 0 $ (resp.\ $ 1$), 
and then $ 2h $ (resp.\ $ 2h - 1$) is the same map,

\Ndt - case (ii) fails: the argument only works if the path $ a\c \cI \to X $ can be splitted as $ a = a' * a''$.

\LB

\Ndt (c) Finally, we have proved that: {\em if} $ X $ {\em is a flexible c-space}, the functor 
$ X \ti - $ (or equivalently the functor $ - \ti X$) preserves the pushout \eqref{4.2.1}, in $ \cTop$.

	We shall use the fact that the product by $ \uI $ works, using flexible homotopies of general 
paths. (In this case the preservation of all colimits already follows from $ \uI $ being 
exponentiable in $ \cTop$, by Theorem \ref{3.9}(b).)

\subsection{Double paths and 2-paths}\label{4.3}
(a) A map $ H\c \cI \ti \uI \to X $ represents a flexible homotopy between two general paths 
$ H(1 \ti \dd^\al)\c \cI \ti \sing \to X$
    \begin{equation}
H\c a \to b\c \cI \to X,   \q\;\;   H(-, 0)  =  a,   \;\;\;   H(-, 1)  =  b.
    \label{4.3.1} \end{equation}

	$H $ will be called a {\em hybrid double path}, as it is parametrised on the hybrid square 
$ \cI \ti \uI$; the latter is less fine than $ \cI^2$, and gives a stronger condition on $ H$.

\Ndt (b) In particular, a {\em hybrid 2-path} is a map $ H\c \cI \ti \uI \to X $ which is constant on the vertical 
faces of the square; it will be written as $ H\c a \to_2 b\c x \to y$, where the c-paths $ a $ and $ b $ 
are its horizontal faces: see figure \eqref{3.4.4}.

\Ndt (c) Marginally, we also consider {\em general double paths} $ H\c \cI^2 \to X$; such a map gives 
a {\em general 2-path} $\, H\c a \to_2 b\c x \to y \,$ if the vertical faces $H(\al, -)$ are constant.

\Ndt (d) All these notions extend the corresponding ones for $\dTop$: if $ X $ is a d-space

\LL

\ndt - a general path $ \cI \to X $ is the same as a d-path $ \uI \to X $ (since $ (\cI)\hat{\,} = \uI$),

\Ndt - a hybrid 2-path $ \cI \ti \uI \to X $ is the same as a 2-path $ \uI^2 \to X $ 
(because $ (\cI \ti \uI)\hat{\,} = \uI^2$, by \ref{2.6}).

\LB

\Ndt (e) A homotopy $ \ph\c f \to g\c X \to Y $ will be said to be {\em strict} if it is constant at each flexible 
point $ x \in X$: the path $ \ph(x, -)\c f(x) \to g(x) $ is a trivial loop in $ Y$. Then $ f $ and $ g $ have the 
same restriction $ f_0 = g_0\c |X|_0 \to |Y|_0 $ to the flexible supports.

	This notion is of interest if the c-space $ X $ has few flexible points, while it gives a trivial 
homotopy if $ X $ is a d-space. In particular, a strict (resp.\ strict flexible) homotopy between c-paths 
is the same as a 2-path (resp.\ hybrid 2-path).

	Strict homotopies will be important for border flexible c-spaces, in Part III.

\subsection{Complements}\label{4.4}
Let $ H\c \cI \ti \uI \to X $ be a hybrid 2-path between c-paths $ a, b\c x \to y$.

\Ndt (a) At any $ t \in \uI $ (always a flexible point) we get an {\em intermediate} c-path
    \begin{equation}
H_t  =  H(-, t)\c \cI \to X,   \qq   H_t\c x \to y,
    \label{4.4.1} \end{equation}
that varies continuously (in the path space $X^\bbI$) from $ H_0 = a $ to $ H_1 = b$. These paths cover 
$\Im H$, in the sense that $\Im H = \Cup_t \, \Im H_t$. This proves that a hybrid 2-path $ H $ between c-paths 
$ x \to y $ satisfies
    \begin{equation}
\Im H  \sub  |X(x, y)|,
    \label{4.4.2} \end{equation}
where $ X(x, y) $ denotes the c-subspace of $ X $ formed by the union of the images of all c-paths 
$ x \to y $ in $ X$.

\Ndt (b) On the other hand, the only flexible points of $ \cI $ are the endpoints, where we get 
$ e_x $ and $ e_y$. At any other $ s \in \cI $ we get a topological path $ H(s, -)\c$ $a(s) \to b(s)$, 
which is actually increasing for the {\em extended path preorder} $ x \precc x' $ of $ X $ 
(determined by the existence of a generated path $ x \to x' $ in $ \hat{X}$, see I.1.7(c)).

	In fact the hybrid double path $ H $ `is' also a d-map $ (\cI \ti \uI)\hat{\,} \to \hat{X}$, and 
$ (\cI \ti \uI)\hat{\,} = \uI^2 $ (by \ref{2.6}), so that $ H(s, -) $ `is' a d-map $ \uI \to \hat{X}$.

\Ndt (c) Let $ x, y $ be two flexible points of the c-space $ X$. The subspace $ X(x, y) $ covered 
by the c-paths $ x \to y $ of $ X $ gives a functor $ \uPi_1(X(x, y)) \to \uPi_1(X) $ which is 
bijective for the arrows $ x \to y$.

	Surjectivity is obvious. Injectivity is a trivial consequence of \eqref{4.4.2}: every hybrid 
2-path between c-paths $ x \to y $ has image in $ X(x, y)$.

\Ndt (d) An easy construction shows that a general 2-path $ H\c \cI^2 \to X $ between c-paths 
$ x \to y $ also has a continuous family of intermediate paths $ Hw_t\c \cI \to X$, which vary 
from $ Hw_0 = a * e_y $ to $ Hw_1 = e_x * b$, and cover $\Im H$. Therefore, property 
\eqref{4.4.2} is still satisfied.

\subsection{The concatenation of hybrid double paths}\label{4.5}
Hybrid double paths have a horizontal concatenation $ H *_1 K$, and no vertical concatenation.

\Ndt (a) {\em Horizontal concatenation}. The concatenation pushout \eqref{4.2.1} is preserved by the 
product by $ \uI$, a flexible c-space, as we have seen in \ref{4.2}(c)
%
    \begin{equation} \begin{array}{c} 
    \xymatrix  @C=10pt @R=8pt
{
~\cI~   \ar[rr]^-{\ddp \ti \, 1}   \ar[dd]_(.4){\ddm \ti \, 1}    &&   ~\cI \ti \uI~    \ar[dd]^(.4){c^- \ti 1}
\\ 
&&   \ar@{.}@/_/[ld] 
\\ 
~\cI \ti \uI~   \ar[rr]_(.45){c^+ \ti 1}    &&   ~\cJ \ti \uI~
}
    \label{4.5.1} \end{array} \end{equation}

	Given two horizontally consecutive hybrid double paths $H, K\c \cI \ti \uI$ $ \to X$, this pushout 
gives a map $ L\c \cJ \ti \uI \to X$. Using the reshaping $ \ka' = \ka \ti \ul{i}\c \cI \ti \uI \to \cJ \ti \uI$, 
we define the horizontal pasting of $ H $ and $ K $ as the composed map
    \begin{equation}
H *_1 K  =  L\ka'\c \cI \ti \uI \to X.
    \label{4.5.2} \end{equation}

\Ndt (b) There is no vertical concatenation of hybrid double paths: this is proved in \ref{4.7}(b).

\Ndt (c) The existence of a hybrid 2-path $ a \to_2 b $ in $ X $ is a reflexive relation (between paths in 
$ X$, with the same endpoints), consistent with path concatenation by the horizontal concatenation (a).

\vskip 2pt

	Transitivity fails, as stated in (b), but this can be overcome: we write as $ a' \precc_2 a'' $ the 
preordering spanned by the previous relation: there exists a finite sequence $ a' \to_2 a_1 \to_2 a_2 \, ... 
\to_2 a'' $ of hybrid 2-paths. We write as $ a' \eqq a'' $ the equivalence relation, called {\em flexible 
2-equivalence}, spanned by the latter (or by $ \to_2$): there exists a finite sequence 
$$
a' \;\;  \longrightarrow_2  \;\;  a_1 \;\;_2\!\!\longleftarrow  \;\;  a_2  \;\; ...  \;\;   \longrightarrow_2  \;\;   a''
$$
of hybrid 2-paths, forward or backward. Both these relations  are consistent with concatenation.

	A class of paths $ [a] $ {\em up to 2-equivalence} will be a class of this equivalence relation, and 
an arrow of the fundamental category of $ X$, in \ref{5.1}.

\Ndt (d) The failure of transitivity is not a real problem here: in any case we must use the equivalence 
relation generated by the relation $ a' \to_2 a''$, and it makes little difference whether the latter is 
transitive or not. (It will make a difference in the general theory of homotopies, where we do not want 
to miss direction: see Part III.)
	
	On the other hand, the failure of point (a) would have been crucial: the congruence of categories 
generated by a relation between parallel arrows which is not consistent with composition involves all 
the objects, and is too `uncontrolled' to allow non-trivial calculations.

\begin{small}

\Ndt (e) If we only work on $ \cI$, by general homotopies of general paths, horizontal and vertical 
concatenation both fail (they are equivalent, by symmetry): general double paths, parametrised on 
$\cI^2$, are not closed under concatenation. This is proved in \ref{4.7}(c).

	On the other hand, if we only work on $ \uI$, by flexible homotopies of flexible paths, we just get 
the fundamental category $ \uPi_1(\Fl X) $ of the flexible part of $ X$: in fact a map $ \uI \to X $ is a 
d-path in $ \Fl X$, and a map $ \uI^2 \to X $ is a double d-path in $ \Fl X$.

\end{small}

\subsection{Lemma}\label{4.6}
{\em
Let $ a\c x \to y $ be a path in the c-space $ X$.

\Ndt (a) If $ \rho, \si\c \bbI \to \bbI $ are global reparametrisations and $ \rho \le \si$, there exists a 
hybrid 2-path:
    \begin{equation}
a\rho \, \to_2 \, a\si.
    \label{4.6.1} \end{equation}

\Ndt (b) In particular, using the global reparametrisations $ \rho, \si\c \bbI \to \bbI $ of \eqref{1.5.1} and 
\eqref{1.5.2}, we get hybrid 2-paths of {\em acceleration} and {\em associativity}:
    \begin{equation}
e_x * a \, \to_2 \, a,   \qq\;\;   a \, \to_2 \, a * e_y,
    \label{4.6.2} \end{equation}
    \begin{equation}
a * (b * c) \, \to_2 \, a * b *c,   \q   a * b * c \, \to_2 \, (a * b) * c.
    \label{4.6.3} \end{equation}

\Ndt (c) For any global reparametrisation $ \rho\c \bbI \to \bbI $ there exist hybrid 2-paths
    \begin{equation}
a\rho' \, \to_2 \, a \, \to_2 \, a\rho'',     \q\;\;     a\rho' \, \to_2 \, a\rho \, \to_2 \, a\rho''
    \label{4.6.4} \end{equation}
proving that $ a \, \eqq a\rho$.
}
\begin{proof}
It is a consequence of Lemma \ref{3.7}, once we verify that the increasing interpolation 
$\Phi\c \rho \to_2 \si$ constructed there
    \begin{equation}
\Phi\c \bbI^2 \to \bbI,   \q\;\;     \Phi(s, t)  =  (1 - t)\rho(s) + t\si(s),
    \label{4.6.5} \end{equation}
is a map $ \cI \ti \uI \to \cI$. 

	In fact, the hybrid square $ \cI \ti \uI$ has three kind of controlled paths 
$ \cI \to \cI \ti \uI $ (see I.2.7(c)), of the form $ \lan 0, v\ran $, or $ \lan 1, v\ran $, or $ \lan u, v\ran $, 
with increasing functions $ u, v $ and $ u $ surjective
%
    \begin{equation} 
\xy <.5mm, 0mm>:
%
(-20,22) *{}; (0,-22) *{};
(50,-10) *{\cI \ti \uI};
@i@={(0,-15), (30,-15), (30,15), (0,15)},
s0="prev"  @@{;"prev";**@{.}="prev"};
(0, -10); (30, 10) **\crv{(5,5)&(20,0)},
\POS(15,2.5) \are+(4,1.2), 
\POS(0,0) \ar+(0,10), \POS(30,-10) \ar+(0,10),
\endxy
    \label{4.6.6} \end{equation}

	Now the functions $ \Phi(0, v(t)) $ and $ \Phi(1, v(t)) $ are constant at $ 0 $ and $ 1$, respectively, 
while $ \Phi(u(t), v(t)) $ is increasing from 0  to 1.	
\end{proof}
%

\subsection{Concatenating double paths}\label{4.7}
The concatenation of double paths in $\dTop$ was recalled in \ref{3.4}. We end this section by 
analysing what happens in $\cTop$ for hybrid, general and flexible double paths. The particular case 
of 2-paths is dealt with in \ref{4.8}.

\Ndt (a) ({\em Horizontal concatenation of hybrid double paths}) This case has already been discussed 
in \ref{4.5}(a): we have a concatenation pushout \eqref{4.5.1}, which allows us to define the horizontal 
concatenation $ H *_1 K $ of horizontally consecutive {\em hybrid double paths}. 

\Ndt (b) ({\em Vertical concatenation of hybrid double paths}) In the following pushout, $P$ is a 
c-structure on $\bbI^2$ generated by the structural maps $1 \ti c^\pm$ 
%
    \begin{equation} \begin{array}{c} 
    \xymatrix  @C=7pt @R=8pt
{
~\cI~   \ar[rr]^-{1 \times \ddp}   \ar[dd]_(.4){1 \times \ddm}    &&   ~\cI \ti \uI~    \ar[dd]^(.4){1 \times c^-}
\\ 
&&   \ar@{.}@/_/[ld] 
\\ 
~\cI \ti \uI~   \ar[rr]_-{1 \times c^+}    &&   ~P~
}
    \label{4.7.1} \end{array} \end{equation}

	This structure is not $\cI \ti \uI$ (in fact, it is strictly finer): we prove below, in point (e), that the 
diagonal of the square is not controlled in $P$. Therefore:

\LL

\ndt - the interval $ \cI $ is not exponentiable in $ \cTop$, since the functor $ \cI \ti - $ does not 
preserve the concatenation pushout \eqref{3.2.1} (of the d-interval),

\Ndt - hybrid double paths are not closed under vertical concatenation (within topological double 
paths), and flexible homotopies cannot be concatenated in $ \cTop$.

\LB

\begin{small}

\skp	Indeed, the structural maps $\cI \ti \uI \to P$ are vertically consecutive hybrid double paths, 
and flexible homotopies (of general paths). Their `topological' vertical concatenation (in $ \Top$) is 
$\id P$, which is not a c-map $\cI \ti \uI \to P$.

\end{small}

\Ndt (c) ({\em Concatenation of general double paths}) By symmetry it is sufficient to consider the 
vertical case. The following pushout gives a c-structure $Q$ on the square, different from $\cI^2$
%
    \begin{equation} \begin{array}{c} 
    \xymatrix  @C=12pt @R=8pt
{
~\cI~   \ar[rr]^-{1 \times \ddp}   \ar[dd]_(.43){1 \times \ddm}    &&   ~\cI^2~    \ar[dd]^(.43){1 \times c^-}
\\ 
&&   \ar@{.}@/_/[ld] 
\\ 
~\cI^2~   \ar[rr]_-{1 \times c^+}    &&   ~Q~
}
    \label{4.7.2} \end{array} \end{equation}

	Again, the diagonal of the square is not controlled in $Q$, which is finer than 
$P$ in pushout \eqref{4.7.1}, as $\cI^2$, in the present pushout, is finer than $\cI \ti \uI$.  

\begin{small}
\vskip 3pt

	We also note that the concatenation pushout $ \cJ$, in \eqref{4.2.1}, is not preserved by the 
cylinder functor $ - \ti \cI $: in fact, $\cJ \ti \cI $ is distinct from $Q$, by the same reason.

\end{small}
\vskip 2pt

	As a consequence:
\LL

\ndt - general double paths are not closed under vertical concatenation (within topological double 
paths),

\Ndt - general homotopies cannot be concatenated in $ \cTop$.

\LB
\vskip 3pt

\Ndt (d) ({\em Concatenation of flexible double paths}) Flexible squares have a horizontal and vertical 
concatenation in $ \cTop $. Indeed a c-map $\uI^2 \to X $ is the same as a d-map $ \uI^2 \to \Fl X$, and 
we are just considering d-squares in the flexible parts.

\Ndt (e) Finally we verify that, in the c-space $P$ of pushout \eqref{4.7.1}, the diagonal map $d$ of 
$\bbI^2$ is not a c-path
%
    \begin{equation} 
\xy <.5mm, 0mm>:
%
(-20,22) *{}; (0,-22) *{};
(-3,-4) *{\sst{d}}; (22,-10) *{\bbI^2}; 
(75,0) *{\sst{p'}}; (116,0) *{\sst{p''}}; (85,11) *{\sst{a}} ; (90,-10) *{\sst{b}} ; (89,3) *{\sst{c}}; 
(122,-10) *{P};
@i@={(-20,-15), (10,-15), (10,15), (-20,15)},
s0="prev"  @@{;"prev";**@{.}="prev"};
@i@={(80,-15), (110,-15), (110,15), (80,15)},
s0="prev"  @@{;"prev";**@{.}="prev"};
(80, -15); (110, 0) **\crv{(85,0)&(100,-10)}, (80, 0); (110, 15) **\crv{(85,15)&(100,5)},
\POS(95,-5.3) \are+(4,.2), \POS(95,9.7) \are+(4,.2), 
\POS(-20,-15) \arl +(30,30), \POS(80,-15) \arl+(0,15), 
\POS(80,0) \arl+(30,0), \POS(110,0) \arl+(0,15), 
\POS(-5,0) \ar+(3,3), \POS(95,0) \are+(4,0), 
\endxy
    \label{4.7.3} \end{equation}

	In fact any c-path $(0, 0) \to (1, 1)$ in $P$, being generated by the c-paths of the lower and upper 
half, must go through the point $ p' = (0, 1/2) $ (as $ a $ in the right figure above), or $ p'' = (1, 1/2) $ 
(as $ b$), or both (as $ c$).

\subsection{General 2-paths}\label{4.8}
We have seen that general double paths are not closed under horizontal or vertical concatenation 
(within topological double paths), but the argument is not really conclusive for our goal: we have to 
show that general {\em 2-paths} are not closed under horizontal concatenation.

	We start from the following pushout (symmetric with respect to \eqref{4.7.2})
%
    \begin{equation} \begin{array}{c} 
    \xymatrix  @C=12pt @R=8pt
{
~\cI~   \ar[rr]^-{ \ddp \times 1}   \ar[dd]_(.43){ \ddm \times 1}    &&   ~\cI^2~    \ar[dd]^(.43){c^- \times 1}
\\ 
&&   \ar@{.}@/_/[ld] 
\\ 
~\cI^2~   \ar[rr]_-{c^+ \times 1}    &&   ~Q'~
}
    \label{4.8.1} \end{array} \end{equation}

	We take the quotient $ Q'/R$, modulo the equivalence relation that collapses each of the thick 
segments to a point, in the left figure below
%
    
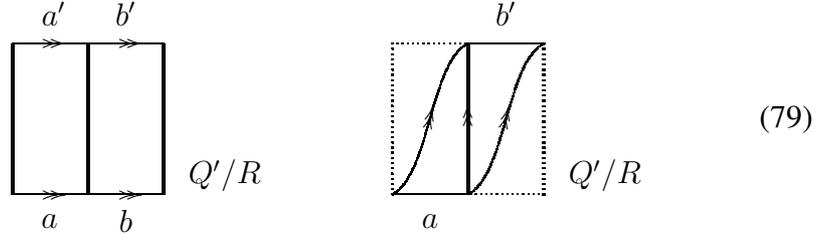
\begin{equation} 
\xy <.5mm, 0mm>:
%
(0,30) *{}; (0,-30) *{};
(0,-27) *{a}; (20,-27) *{b}; (1,28) *{a'}; (20,28) *{b'}; (46,-15) *{Q'/R}; 
(100,-27) *{a}; (120,28) *{b'}; (146,-15) *{Q'/R};
@i@={(-10,-20), (30,-20), (30,20), (-10,20)},
s0="prev"  @@{;"prev";**@{-}="prev"}; 
@i@={(90,-20), (130,-20), (130,20), (90,20)},
s0="prev"  @@{;"prev";**@{.}="prev"}; 
(90, -20); (110, 20) **\crv{(100,-15)&(100,15)}, (110, -20); (130, 20) **\crv{(120,-15)&(120,15)}, 
\POS(90,-20) \arl+(20,0), \POS(110,20) \arl+(20,0), 
\POS(0,-20) \are+(4,0), \POS(20,-20) \are+(4,0), 
\POS(0,20) \are+(4,0), \POS(20,20) \are+(4,0), 
\POS(100,0) \are+(.7,3), \POS(120,0) \are+(.7,3), \POS(110,0) \are+(0,4), 
\POS(-9.6,-20) \arl+(0,40), \POS(9.8,-20) \arl+(0,40), 
\POS(10.2,-20) \arl+(0,40), \POS(29.6,-20) \arl+(0,40),
\POS(90.3,-20) \arlp+(0,40), \POS(109.8,-20) \arl+(0,40), 
\POS(110.2,-20) \arl+(0,40), \POS(129.7,-20) \arlp+(0,40),
\endxy
    \label{4.8.2} \end{equation}

	Now $c^- \ti 1$ and $c^+ \ti 1$ induce general 2-paths $ a \to_2 a'$ and $ b \to_2 b'$, respectively. 
Topologically, their horizontal concatenation is the canonical projection $ p\c Q' \to Q'/R$, and this is not 
a c-map $\cI^2 \to Q'/R$, by an argument similar to that of \eqref{4.7.3}.

\begin{small}
\vskip 3pt

	The diagonal $ d\c \cI \to \cI^2 $ is projected to a path $(0, 0) \to (1, 1)$ in the support of 
$ Q'/R$, which is not a c-path: to qualify as such it should admit as a restriction either $ a $ or $ b' $ or 
both, as shown in the right figure above.

\end{small}

\subsection{Trivial loops}\label{4.9}
Finally we want to make clear a point concerning the definition of c-spaces in \ref{1.1}: if we replace 
axiom (csp.0) about constant paths with (dsp.0) (all trivial loops are selected), we get a structure -- 
let us say a {\em \ul{c}-space} -- where our construction of the fundamental category fails, by the failure of the 
horizontal concatenation of hybrid double paths, in \ref{4.7}(a).

	In fact, the standard \ul{c}-interval $\ul{\rc}\bbI$ is the euclidean interval $ [0, 1] $ with the new 
structure generated by the identity mapping $ \ul{i}$: the selected paths are the surjective increasing 
maps $ \bbI \to \bbI $ and all the constant ones.

	The new hybrid square has two kind of c-paths: the increasing vertical paths 
$t \mapsto (s_0, v(t))$ and all increasing paths with surjective first projections.

	Now the double path $ \Phi\c \bbI^2 \to \bbI $ defined in \eqref{4.6.5} is not a map 
$\ul{\rc}\bbI \ti \uI \to \ul{\rc}\bbI$: on a vertical path $ t \mapsto (s_0, t) $ in $\ul{\rc}\bbI \ti \uI $ we get a 
path $ \Phi(s_0, -)\c \bbI \to \bbI $ from $ \rho(s_0) $ to $ \si(s_0)$, which -- generally -- is not selected 
in $\ul{\rc}\bbI$: it is neither surjective nor constant.

\section{The fundamental category of controlled spaces}\label{s5}

		We introduce the fundamental category $\uPi_1(X)$ of a c-space, extending the fundamental 
category of a d-space -- and therefore the fundamental groupoid $\Pi_1(-)$ of a topological space.

	There are also fundamental categories $ \uPi_1(\Fl X) $ and $ \uPi_1(\hat{X}) $ induced by the 
coreflector and the reflector of d-spaces, and linked to the previous one by obvious natural transformations. 
But these `induced' functors miss the critical aspects of c-spaces.

	The invariance of $\uPi_1$ up to flexible homotopies, proved in Theorem \ref{5.4}, will be linked 
to directed homotopy equivalence of categories \cite{G3}, in Part III.

	After Theorem \ref{5.3} we already have some of the main ingredients to compute the fundamental 
category of a c-space; one can read at that point the elementary computations of $\uPi_1(X)$ in \ref{5.9} 
(which only marginally rely on Theorem {5.8}).

\subsection{Definition {\rm (The fundamental category of a c-space)}}\label{5.1}
The {\em fundamental category} $ \uPi_1(X) $ of a controlled space consists of the following items:

\LL

\ndt - the vertices are the flexible points of $ X$,

\Ndt - the arrows are the equivalence classes $\, [a]\c x \to y \,$ of general paths $ a\c x \to y$, up to the 
2-equivalence relation $a' \, \eqq \, a''$ spanned by the hybrid 2-paths (see \ref{4.5}(c)),

\Ndt - the composition is induced by concatenation of general paths, and the identities are induced by 
trivial loops
    \begin{equation}
[a][b]  =  [a * b],   \q\;\;\;   1_x  =  [e_x].
    \label{5.1.1} \end{equation}

\LB
\begin{small}

	The operation is well defined because 2-equivalence is consistent with path concatenation. 
Identities and associativity follow from hybrid 2-paths of acceleration and associativity, in Lemma 
\ref{4.6}.

\end{small}
\vskip 2pt

	A d-map $ f\c X \to Y $ produces a functor between small categories
    \begin{equation} \begin{array}{c}
f_*  =  \uPi_1(f)\c \uPi_1(X) \to \uPi_1(Y),
\\[5pt]
f_*(x)  =  f(x),   \q   f_*[a]  =  [fa]   \qq   (x \in |X|_0),
    \label{5.1.2} \end{array} \end{equation}
and we have a functor
    \begin{equation}
\uPi_1\c \cTop \to \Cat,
    \label{5.1.3} \end{equation}
which extends $ \uPi_1\c \dTop \to \Cat$, by \ref{4.3}(d).

	At a flexible point $ x_0 $ of the c-space $ X$, the endoarrows of $ \uPi_1(X) $ form the 
{\em fundamental monoid}
    \begin{equation}
\upi_1(X, x_0)  =  \uPi_1(X)(x_0, x_0).
    \label{5.1.4} \end{equation}

	Extending a definition of d-spaces (\cite{G3}, 3.2.7), we say that a c-space $ X $ is {\em 1-simple} 
if its fundamental category is a preorder (i.e.\ the category associated to a preordered set). This means 
that every hom-set $ \uPi_1(X)(x, x') $ has at most one arrow: there is one if there is a c-path $ x \to x' $ 
in $ X$, and none otherwise. We shall see that many basic c-spaces are of this kind.

\subsection{Induced fundamental categories}\label{5.2}
We can also use in $ \cTop $ {\em the homotopy theory of} $\dTop$, through the reflector 
$(\hat{\;\;})\c \cTop \to \dTop $ and the coreflector $ \Fl\c \cTop \to \dTop$. 

	We obtain thus two `induced' functors, the {\em fundamental category of generated paths} 
$ \uPi_1(\hat{X})$
    \begin{equation}
\uPi_1(\hat{\;\;})\c \cTop \to \Cat,   \q\;\;\;   X \mapsto \uPi_1(\hat{X}),  \;\;
    \label{5.2.1} \end{equation}
and the {\em fundamental category of flexible paths} $\uPi_1(\Fl X)$
    \begin{equation}
\uPi_1\Fl\c \cTop \to \Cat,   \qq   X \mapsto \uPi_1(\Fl X).
    \label{5.2.2} \end{equation}

	We recall that the support of $ \Fl X $ is the subspace $ |X|_0 $ of flexible points, which are the 
vertices of the categories $\uPi_1(\Fl X)$ and $\uPi_1(X)$.

	The canonical embeddings $ \Fl X \to X \to \hat{X} $ (the counit and unit of the adjunctions) give 
two natural transformations, with components
    \begin{equation}
\uPi_1(\Fl X) \, \longrightarrow \, \uPi_1(X)  \, \longrightarrow \, \uPi_1(\hat{X}).
    \label{5.2.3} \end{equation}

	If $ X $ is flexible we get two identities. These functors need not be faithful (see Part III), but 
Theorem \ref{5.3}(b) shows an important case where the second is a full embedding.

	As we shall see in \ref{5.9}, these three categories are strongly different on rigid objects. Both the 
`induced' functors $ \uPi_1(\Fl X) $ and $ \uPi_1(\hat{X}) $ are unable to analyse the critical features 
of c-spaces: the former leaves out the non-flexible paths and the latter makes all paths flexible.

	For every c-space $X$, the underlying topological space $ |X| $ has a natural c-structure, and 
the reshaping $X \to |X|$ gives another natural transformation $ \uPi_1(X) \to \Pi_1(|X|) $ with values in 
the ordinary fundamental category of the support. There is thus a commutative diagram
    \begin{equation} \begin{array}{c} 
    \xymatrix  @C=20pt @R=20pt
{
~\uPi_1(\Fl X)~   \ar[r]   \ar[d]    &   ~\uPi_1(X)~   \ar[r]   \ar[d]    &   ~\uPi_1(\hat{X})~    \ar[d]
\\ 
~\Pi_1(|X|_0)~   \ar[r]     &   ~\Pi_1(|X|)~   \aru[r]    &   ~\Pi_1(|X|)~
}
    \label{5.2.4} \end{array} \end{equation}
%
%

\subsection{Theorem {\rm (Weak flexibility and fundamental category)}}\label{5.3}
{\em
(a) If $ X $ is full in the c-space $ X'$ (see  \ref{2.1}(a)), $\uPi_1(X)$ is the full subcategory 
of $\uPi_1(X') $ with vertices in $ |X|_0$.

\Ndt (b) If the c-space $ X $ is preflexible (i.e.\ full in the generated d-space $ \hat{X}$, see \ref{2.1}(b)) 
the component $ \uPi_1(X) \to \uPi_1(\hat{X}) $ is a full embedding. In other words $ \uPi_1(X) $ 
is the full subcategory of $\uPi_1(\hat{X})$ with vertices in $ |X|_0$.
}
\begin{proof}
It is sufficient to prove (a). Let $ x, y $ be two flexible points of the c-space $ X$. The subspace $ X(x, y) $ 
covered by the c-paths $ x \to y $ of $ X $ coincides with $ X'(x, y)$, and the component 
$\uPi_1(X'(x, y)) \to \uPi_1(X')$ is bijective for the arrows $x \to y$, by \ref{4.4}(c). 
\end{proof}
%

\subsection{Theorem {\rm (Homotopy invariance, II)}}\label{5.4}
{\em
(a) The fundamental category of c-spaces is invariant in the following sense: a {\em flexible} homotopy 
$ \ph\c f \to g\c X \to Y $ of c-spaces induces a natural transformation
    \begin{equation} \begin{array}{c}
\ph_*\c f_* \to g_*\c \uPi_1(X) \to \uPi_1(Y),
\\[5pt]
\ph_*(x)  =  [\ph_x]\c f(x) \to g(x),   \q   \ph_x  =  \hat{\ph}(x \ti \ul{i})\c \sing \ti \uI \to Y,
    \label{5.4.1} \end{array} \end{equation}
where, for every $x \in |X|_0$, $\, \ph_x$ is a flexible path in $Y$, from $ f(x) $ to $ g(x)$.

\Ndt (b) If $ \ph\c f \to g $ is a strict flexible homotopy (see \ref{4.3}(e)), $ \ph_* $ is the identity of 
$ f_* = g_*$.

\Ndt (c) A whisker composition $\psi = k \ci \ph \ci h\c kfh \to kgh $ of flexible homotopies and maps is 
taken to the corresponding whisker composition of natural transformations and functors.
}
\begin{proof}
(a) Naturality works as in Theorem \ref{3.6}(a), for d-spaces. We have only to check that, in the present 
setting, the three double paths used in \eqref{3.6.5} are maps $ \cI \ti \uI \to Y$, that is {\em flexible} 
double paths. (Of course this would not work for a general homotopy.)

	First, the double path $ \ph_a = \ph(a \ti \ul{i}) $ is a composite $ \cI \ti \uI \to X \ti \uI \to Y$. 
Secondly, the double path $ \ph_x g^+ $ is a composite $ \uI^2 \to \uI \to Y$, because $ \ph_x $ is a 
flexible path; but $ \cI \ti \uI $ is finer than $ \uI^2 $ and we are done. The third double path 
$ \ph_{x'}g^- $ works in the same way.

\Ndt (b) An obvious consequence. (c) By the same argument as in \ref{3.6}(b).
\end{proof}
%

\subsection{Theorem}\label{5.5}
{\em
(a) If $ a $ is a flexibly reversible c-path in the c-space $ X $ (see \ref{4.1}), the arrows $ [a] $ and 
$ [ar] $ are inverse to each other in $ \uPi_1(X)$.

	This need not be true if $a$ is merely flexible and reversible: see Part III.

\Ndt (b) If $ \ph\c f \to g\c X \to Y $ is a flexibly reversible homotopy of c-spaces (see \ref{4.1}), the 
natural transformation $ \ph_*\c f_* \to g_* $ is invertible.

	More generally, the same holds if $ \ph $ is a flexible homotopy of c-spaces and, for every 
$ x \in X$, the c-path $ \ph_x = \hat{\ph}(x \ti \ul{i})\c \sing \ti \uI \to Y $ is flexibly reversible in $ Y$.
}
\begin{proof}
The path $ a\c x \to x' $ is supposed to be a map $ a\c \bbI^\sim \to X$. The double paths 
$H(s, t)  =  a(s \me t)$ and $K(s, t)  =  a((1 - s) \me t)$
%
    \begin{equation} 
\xy <.45mm, 0mm>:
%
(0,23) *{}; (140,-23) *{};
(25,0) *{H}; (3,0) *{\sst{e_x}}; (43,0) *{\sst{a}}; (55,0) *{K}; (77,0) *{\sst{e_x}}; 
(25,-22) *{\sst{e_x}}; (55,-22) *{\sst{e_x}}; (25,20) *{\sst{a}}; (54,20) *{\sst{ar}}; 
(120,10) *{H  =  ag^+,}; 
(132,-10) *{K  =  ag^+ (r \ti 1),};   
@i@={(10,-15), (70,-15), (70,15), (10,15)},
s0="prev"  @@{;"prev";**@{-}="prev"};
\POS(10.4,-15) \arl+(0,30), \POS(10,-14.5) \arl+(60,0), \POS(69.6,-15) \arl+(0,30), 
\POS(40,-15) \arl+(0,30), 
\endxy
    \label{5.5.1} \end{equation}
are maps $ \bbI^{\sim 2} \to \bbI^\sim \to X$, and therefore $ \cI \ti \uI \to X$, that is hybrid double paths.

	The latter have a horizontal concatenation $ H *_1 K$, which is a hybrid 2-path and proves that 
$(a * ar) \eqq e_x$. Applying the same result to the path $ar\c x' \to x$, we get 
$(ar * a) \eqq e_{x'}$.

\Ndt (b) An obvious consequence of Theorem \ref{5.4} and the previous point: in the given hypotheses 
each component $ \ph_*(x) = [\ph_x]\c f(x) \to g(x) $ of the natural transformation $ \ph_* $ is invertible.
\end{proof}
%

\subsection{Theorem {\rm (Products and sums)}}\label{5.6}
{\em
(a) The functor $ \uPi_1\c \cTop \to \Cat $ preserves arbitrary products and sums.

\Ndt (b) A product or sum of a family $(X_i)$ of 1-simple c-spaces is 1-simple. Conversely, if their product 
or sum is 1-simple, all $ X_i $ are also -- provided there is no empty factor in the case of the product. 
}
\begin{proof}
(a) One can apply the same argument that works in $ \Top $ (and $\dTop$). In a cartesian product of 
c-spaces, paths and hybrid double paths are detected by the cartesian projections. In a sum, they live 
in one `component'.

\Ndt (b) A consequence of (a). In fact, property (b) holds in $ \Cat$, interpreting `1-simple' as being a 
preorder. Moreover, preorders form a reflective and coreflective subcategory of $\Cat$, so that limits 
and colimits agree.

\begin{small}

	(In the elementary homotopy theory of categories described in \cite{G3}, 1.1.6, 
the fundamental category of a category is the category itself.)

\end{small}
\end{proof}
%

\subsection{Covering maps}\label{5.7}
(a) The theory of covering maps for topological spaces, a classical topic of Homotopy Theory, can be 
found in many books on Algebraic Topology, like \cite{Ha, HiW, Hu, Sp}.

	We recall that a map $ p\c X \to Y $ of topological spaces is called a {\em covering map}, or a 
covering projection, if every point $ y \in Y $ has an open neighbourhood which is `evenly covered' by 
$ p$: this means that the preimage $ p^{-1}(U) $ is the disjoint union of open subsets of $ X$, each of 
which is mapped homeomorphically onto $U$ by $p$. $X $ is the {\em total space}, or covering space, 
of $p$; $Y$ is the {\em basis}; $F_y = p^{-1}\{y\} $ is the {\em fibre} at $ y \in Y$ -- a discrete 
subspace of $X$. The map $p$ is a surjective local homeomorphism, and an open map.

	A classical example is the {\em exponential map}, forming the `universal covering' of the circle (in 
the complex plane)
    \begin{equation}
p\c \bbR \to \bbS^1,   \q\;\;\;   p(t)  =  e^{2\pi it}.
    \label{5.7.1} \end{equation}

\Ndt (b) Here we only need the basic results on the homotopy lifting property of a covering map 
$p\c X \to Y$. (For a proof without prerequisites one can see \cite{G4}, Theorem 6.2.9.)

\LL

\ndt (i) For every path $ b\c \bbI \to Y $ in the basis and every point $ x_0 $ in the fibre at $ b(0)$, there 
is a unique lifting $ a\c \bbI \to X $ (i.e.\ $ pa = b$) starting at $ x_0$.

\Ndt (ii) For every 2-path $ K\c b \to_2 b' $ (a map $ K\c \bbI^2 \to Y$) and every point $ x_0 $ in the 
fibre at $b(0) = b'(0)$, there is a unique lifting $ H\c \bbI^2 \to X $ (i.e.\ $ pH = K$) to a 2-path 
$ a \to_2 a' $ of paths starting at $ x_0$. (The latter are liftings of the paths $b, b'$.)

\LB
\skp

	Formally, property (i) is a consequence of (ii); nevertheless, stating it is useful.

\Ndt (c) Now, we say that a map $ p\c X \to Y $ of c-spaces is a {\em covering map of c-spaces} if it is 
a covering map of topological spaces and the path-lifting property (i) holds within c-paths: every lifting 
of a c-path of the basis is a c-path of the total space.

	Then the homotopy-lifting property (ii) automatically holds for any hybrid 2-path $K\c \cI \ti \cJ \to Y$,  
because its topological lifting $H\c \bbI^2 \to |X|$ is a c-map $\cI \ti \cJ \to X$ if and only if every c-path 
$\cI \to \cI \ti \cJ \to Y$ is lifted to a c-path of $X$.  

	We also note that the flexible support of $X$ is the union of the fibres of the flexible points of $Y$. 

	Our prime examples are the exponential c-maps
    \begin{equation} \begin{array}{ccc}
p\c \cR \to \cS^1,  &\;\;\;&   p(t) = e^{2\pi it},
\\[5pt]
 p_n\c \rc_n\bbR \to \rc_n\bbS^1,   &&   p_n(t) = e^{2\pi it}.
    \label{5.7.2} \end{array} \end{equation}

	This will be applied in \ref{5.9} to determine $ \uPi_1(\cS^1) $ and $ \uPi_1(\rc_n\bbS^1)$.

\vskip 3pt
\begin{small}

	The second map can be replaced with $\, p'_n\c \cR \to \rc_n\bbS^1$, $\; p'_n(t) = e^{2n\pi it}$.

\end{small}
\skp

	A {\em covering map of d-spaces} is similarly defined, using d-maps and d-paths.

\subsection{Theorem}\label{5.8}
{\em
Let $ p\c X \to Y $ be a covering map of topological spaces.

\Ndt (a) Let $x_0 \in X$, $y_0 = p(x_0)$ and $y \in Y$. The functor 
$p_*\c \Pi_1 (X) \to \Pi_1 (Y) $ induces a bijection of sets
    \begin{equation}
p_*\c \Sum_{x \in F_y} \, \Pi_1 (X)(x_0, x) \, \to \, \Pi_1 (Y)(y_0, y),
    \label{5.8.1} \end{equation}
defined on the disjoint union of the sets $\, \Pi_1 (X)(x_0, x)$, for $x \in F_y$. 

\vskip 4pt
\begin{small}

	In other words, the functor $ p_* $ is surjective on objects and arrows, and faithful.
The composition in $ \Pi_1 (Y) $ is determined by that of $ \Pi_1 (X)$, as any pair of 
composable arrows in $ \Pi_1 (Y) $ can be lifted to a pair of composable arrows in $ \Pi_1 (X)$.

\end{small}

\Ndt (b) For a covering map $ p\c X \to Y $ of c-spaces (as defined above) we have the same result 
for the functor $ p_*\c \uPi_1(X) \to \uPi_1(Y)$.
}
\begin{proof}
The lifting properties (i), (ii) of \ref{5.8} imply that the mapping \eqref{5.8.1} is surjective and injective. 

	In case (b), where 2-equivalence is {\em generated}  by hybrid 2-paths, one should use again the fact that 
all c-paths can be lifted.
\end{proof}
%

\subsection{Elementary calculations}\label{5.9}

Computing the fundamental category will be studied in Part III, but several results can be simply obtained 
using Theorem \ref{5.3}(b) on preflexible c-spaces (and the fundamental category of the generated 
d-spaces, already computed in \cite{G3}), or Theorem \ref{5.8} on covering maps of c-spaces. In 
particular, many basic c-spaces are 1-simple, in the sense of \ref{5.1}: their fundamental category is a 
preorder; of course, the controlled circle is not.
 
	The symbols $\two, \three, \N, \Z, \R$ denote ordered sets and the associated categories; the ordered 
sets 2, 3 and $D|\Z|$ are discrete. $\bbN$ is the one-object category associated to the additive monoid 
of the natural numbers.

\Ndt (a) By Theorem \ref{5.3}(b), the fundamental categories of the preflexible c-spaces $ \cI$, $ \cJ$, $ \cR $ 
are the following ordered sets:
    \begin{equation}
\uPi_1(\cI)  =  \two,    \q   \uPi_1(\cJ)  =  \three,    \q   \uPi_1(\cR)  =  \Z.
    \label{5.9.1} \end{equation}

	For these preflexible spaces the functors \eqref{5.2.3} become inclusions of ordered sets:
    \begin{equation}
2 \to \two \to [0, 1],    \q   3 \to \three \to [0, 2],    \q   D|\Z| \to \Z \to \R.
    \label{5.9.2} \end{equation}

\Ndt (b) The fundamental category of the directed circle $ \uS^1$, as described in \cite{G3}, 3.2.7(d), is the 
subcategory of the groupoid $\Pi_1 \bbS^1$ formed of the classes of anticlockwise paths (in $ \bbR^2$). 
Each monoid $ \upi_1(\uS^1, x) $ is isomorphic to the additive monoid $\bbN$ of natural numbers.

	Applying Theorem \ref{5.3}(b), the fundamental category of the one-stop circle $\cS^1$ amounts to the 
fundamental monoid at the unique flexible point $x_0$ (the point 1 of the complex plane)
    \begin{equation}
\uPi_1(\cS^1)(x_0, x_0)  =  \upi_1(\uS^1, x_0) =  \bbN .
    \label{5.9.3} \end{equation}

	Without using $ \uPi_1(\uS^1))$, this is also proved by Theorem \ref{5.8}(b) applied to the exponential 
map $ \cR \to \cS^1$. Two c-loops $ a, b $ in $ \cS^1 $ are 2-equivalent if and only if they have the same 
length $2k\pi$ (in radians), if and only if they both turn $ k $ times ($k \ge 0$) around the circle, anticlockwise.

\Ndt (c) More generally, the fundamental category of the preflexible $n$-stop circle $ \rc_n\bbS^1 $ (see 
\eqref{1.4.6}) is the full subcategory $\bc_n$ of the fundamental category of 
$(\rc_n\bbS^1)\hat{\,} = \uS^1 = \uR/\bbZ$ on $ n $ flexible points, the vertices $ [i/n] $ (for $ i = 0, ..., n-1$) 
of an inscribed $n$-gon.

	Again, this result can also be obtained using the covering map of c-spaces 
$p_n\c \rc_n\bbR \to \rc_n\bbS$.

\Ndt (d) Applying Theorem \ref{5.6} on cartesian products, we get the following fundamental categories, 
which are (partially) ordered sets
    \begin{equation} \begin{array}{ccc}
\uPi_1(\cI^n)  =  \two^n,   &&   \uPi_1(\cJ^n)  =  \three^n,
\\[5pt]
 \uPi_1(\cI \ti \cJ)  =   \two \ti \three,   &\qq&   \uPi_1(\cR^n)  =  \Z ^n.
    \label{5.9.4} \end{array} \end{equation}
%
%


\vspace{5mm}
\noindent
    Marco Grandis\\
    Dipartimento di Matematica \\
    Universit\`a di Genova \\
    Via Dodecaneso 35 \\
    16146 - Genova, Italy \\
    grandis@dima.unige.it

\end{document}